\documentclass[10pt,reqno]{amsart}
\usepackage{amssymb,accents,color}
\def\noi{\noindent}

\newtheorem{Thm}{Theorem}[section]
\newtheorem{Def}[Thm]{Definition}
\newtheorem{Lm}[Thm]{Lemma}
\newtheorem{Prop}[Thm]{Proposition}
\newtheorem{Cor}[Thm]{Corollary}
\newtheorem{Rem}[Thm]{Remark}

\newtheorem{conjecture}[Thm]{Conjecture}

\def\cal{\mathcal}
\def\Bbb{\mathbb}
\def\mf{\mathfrak}

\def\<{\langle}
\def\>{\rangle}

\def\a{\alpha}
\def\b{\beta}
\def\d{\delta}

\def\th{\theta} 

\def\l{\lambda}
\def\L{\Lambda}

\def\Re{\Bbb R}
\def\F{\Bbb F}

\def\Z{\Bbb Z}

\def\Qu{\Bbb Q}

\def\H{\cal H}

\def\G{\mf h}

\def\W{\ring W}
\def\w{\ring w}

\def\RR{\ring R}
\def\Q{\ring Q}
\def\h{\ring {\cal H}}
\def\n{\ring {\cal N}}
\def\Gc{\ring {\mf h}}
\def\G{\ring {\mf g}}

\def\bo{\ring {\mf b}}

\begin{document}
\subjclass[2000]{33D52, 33D45, 33D80, 20C08, 17B10}
\title[standard bases and macdonald polynomials]
{Standard bases for affine parabolic modules and nonsymmetric Macdonald polynomials}
\author[]{Bogdan Ion}
\thanks{Department of Mathematics, University of Pittsburgh, Pittsburgh, PA 15260.}
\thanks{E-mail address: {\tt bion@pitt.edu} }
\date{May 8, 2006}
\begin{abstract} We establish a connection between (degenerate) nonsymmetric Macdonald polynomials
and standard bases and dual standard bases of maximal parabolic modules of affine Hecke algebras.
Along the way we prove a (weak) polynomiality result for coefficients of symmetric and nonsymmetric Macdonald polynomials.
\end{abstract}
\maketitle

%%%%%%%%%%%%%%%%%%%%%%%%%%%%%%%%%%%%%%%%%%%%%%%%%%%%%%%%%%%%%%%%%%%%%%%%%%%%%
%%%%%%%%%%%%%%%%%%%%%%%%%%%%%%%%%%%%%%%%%%%%%%%%%%%%%%%%%%%%%%%%%%%%%%%%%%%%%
\thispagestyle{empty}
\section*{Introduction}

The symmetric
Macdonald polynomials $P_\l(q,t)$ are a family of Weyl group
invariant functions depending rationally on  parameters $q$ and
$t=(t_s,t_\ell)$, which are associated to any finite, irreducible root system
$\RR$ and are indexed by the anti--dominant elements of the weight
lattice of $\RR$. Introduced originally for root systems of type $A$
as a common generalization of Hall--Littlewood and Jack symmetric
functions it was quickly realized  that they have
deep properties  essentially rooted in two classical
representation--theoretical contexts: the theory of zonal
spherical functions for real (Gel'fand, Harish--Chandra) and $
p$--adic (Satake, Macdonald, Matsumoto) reductive groups. In a
more recent development \cite{sanderson, ion}, the symmetric Macdonald
polynomials were also connected with the representation theory of
affine Kac--Moody groups.

The nonsymmetric Macdonald polynomials $E_\l(q,t)$  (indexed now by the
full weight lattice) are of more recent vintage. They
were introduced by Heckman, Opdam \cite{opdam} (for $t=q^k$ and $q\to 1$),
Macdonald \cite{masterisque} (for $t=q^k$), Cherednik \cite{c3}
(general case, reduced root systems) and Sahi \cite{sahi} (nonreduced
root systems). They turned out to be a crucial tool in all the recent
developments in the theory of orthogonal polynomials, the related
combinatorics and the representation theory of double affine Hecke
algebras (see, for example \cite{c2}). However,  they do not seem to fit easily 
in a classical representation-theoretical framework, the main obstacle 
being precisely their non-invariance under the Weyl group. The first hint at their representation-theoretical
nature  came from \cite{sanderson} (type $A$), \cite{ion} (general type): 
$E_\l(q,\infty)$ are Demazure
characters of basic representations of affine Kac--Moody groups. Furthermore, 
$E_\l(\infty,\infty)$ are Demazure characters of irreducible representations of simple 
algebraic groups \cite{ion}, and $E_\l(\infty,t)$ are, for specific
values of $t$, matrix coefficients  for unramified principal series representations of simple 
$p$--adic groups \cite{ion2}.

The goal of this paper is add on this list another context in which the nonsymmetric polynomials can be interpreted naturally: the Kazhdan--Lusztig theory.
In the symmetric case the connection is well-known: the limits $P_\l(\infty,t)$ of the symmetric Macdonald polynomials are (via the Satake transform) 
the standard basis of the corresponding spherical Hecke algebra. Our main result gives a similar interpretation for
the same limit of the nonsymmetric polynomials: they form the  standard basis of the maximal parabolic module of the corresponding affine 
Hecke algebra. It should be  noted that the Satake transforms of  Kazhdan--Lusztig bases of spherical Hecke algebras are 
irreducible Weyl characters, a fact which follows almost immediately from the knowledge of the standard bases \cite{kato, lusztig2} (see also Theorem \ref{KLantidominant}). Similar explicit 
formulas for Kazhdan-Lusztig bases of maximal parabolic modules are not immediately obtained from similar information, but seem to require 
 new ideas. Conjecture \ref{conjecture} gives some indication of what is expected in this situation.

In brief, the content of the paper is the following. Section \ref{section1} contains well-known combinatorial properties of affine root systems and their Weyl groups. The main result of Section 2 (stated as Theorem \ref{polynomiality} ) is a polynomiality result for certain normalizations of Macdonald polynomials (the normalization factor $e_\l$ is a product of factors of the form $(1-q^at^b)$ with $a,b$ negative integers ).

\noi {\bf Theorem.} {\it
If the root system $\RR$ is reduced then,

\begin{enumerate}
\item For any weight $\l$, the coefficients of $ e_\l E_\l(q,t)$ are
polynomials in $q^{-1}$, $t_s^{-1}$, $t_\ell^{-1}$ with integer coefficients.

\item For any   anti--dominant weight $\l$, the coefficients of $
e_\l P_\l(q,t)$ are polynomials in $q^{-1}$,
$t_s^{-1}$, $t_\ell^{-1}$ with integer
coefficients.

\end{enumerate}
 If the root system $\RR$ is nonreduced, then

\begin{enumerate}
\item[(3)] For any weight $\l$, the coefficients of $ e_\l E_\l(q,t)$ are
polynomials in $~q^{-1}$, $t_{01}^{-1}$, $t^{-1}$, $t_n^{-\frac{1}{2}}t_{01}^{\pm\frac{1}{2}}$,
$t_n^{-\frac{1}{2}}t_{03}^{\pm\frac{1}{2}}$,
$t^{-\frac{1}{2}}_{01}t_{02}^{\pm\frac{1}{2}}$,
$t^{-\frac{1}{2}}_{01}t_{03}^{\pm\frac{1}{2}}$
 with integer
coefficients.

\end{enumerate}}

In general, the result improves on what was previously known about the nature of the coefficients \cite[Corollary 5.3]{c4} (Laurent polynomials in $q,t$), \cite[Section 2.3]{ion} (polynomials in $t^{-1}$, Laurent polynomials in $q$), but it is still far from being optimal since, as observed in many cases, the normalization factor can be further trimmed down without altering the polynomiality of the coefficients. In fact, in type $A$, stronger results are know for both symmetric and nonsymmetric polynomials \cite{kn,sahi2}. The main technical idea of the proof is to use two affine intertwiners (one dependent of $q$, the other independent) in conjunction. The stronger result in type $A$ was handled similarly taking also advantage of the stability  of the relevant polynomials in that case. 

Section 3 is concerned with the limit $q\to \infty$ of nonsymmetric polynomials. The results have been already used in \cite{ion2} to establish a connection between nonsymmetric Macdonald polynomials and matrix coefficients of unramified principal series for reductive $p$--adic groups. Section 4  recalls the basic (maximal parabolic) Kazhdan-Lusztig theory for affine Hecke algebras (in its muti-parameter version \cite{lusztig3}) and explains the connection with the Cherednik-Macdonald theory. Our main result (stated as Theorem \ref{state2}) is the following

\noi {\bf Theorem.} {\it The basis  $\{ {\widetilde E}_\l(q,t)\}_{\l\in P}$ of the parabolic module of the affine Hecke algebra $\H_X^e$ is 
 invariant under the Kazhdan--Lusztig involution. Moreover, 
\begin{enumerate}
\item $\{ {\widetilde E}_\l(\infty,t)\}_{\l\in P}$ is
the standard basis;

\item  $\{ {\widetilde E}_\l(0,t)\}_{\l\in P}$ is
the dual standard basis.
\end{enumerate}}

In type $A$, the result is essentially contained in \cite[Corollary 5.3]{kn} (see also \cite{knop}).

Finally, Section \ref{section5} is examining the interplay between the case when all parameters approach infinity and the case when all parameters approach zero as well as a new geometric interpretation (Theorem \ref{00}) for the polynomials $E_\l(0,0)$.

\subsection*{Aknowledgements}   It is a pleasure to thank Shrawan Kumar and Peter Littelmann whose insights led to Section \ref{geometric}. This work was started while the author was supported in part by a Rackham Faculty Research Fellowship at the University of Michigan and continued under the support of the NSF grant DMS--0536962. 
%%%%%%%%%%%%%%%%%%%%%%%%%%%%%%%%%%%%%%%%%%%%%%%%%%%%%%%%%%%%%%%%%%%%%%%%%%%%%
%%%%%%%%%%%%%%%%%%%%%%%%%%%%%%%%%%%%%%%%%%%%%%%%%%%%%%%%%%%%%%%%%%%%%%%%%%%%%
\section{Preliminaries}\label{section1}

\subsection{Affine root systems} Let $\RR\subset
\Gc^*$ be a finite, irreducible, not necessarily reduced, root
system of rank $n$, and let $\RR^\vee\subset \Gc$ be the dual root
system. We denote by $\{\a_i\}_{1\leq i\leq n}$ a basis of $\RR$
(whose elements will be called simple roots); the corresponding
elements $\{\a_i^\vee\}_{1\leq i\leq n}$ of $\RR^\vee$ will be
called simple coroots. If the root system is nonreduced, let us
arrange that $\a_n$ is the unique simple root such that $2\a_n$ is
also a root. Throughout the paper a special role will be played by
the root $\th$, which is defined as the {\sl highest short root}
in $\RR$ if the root system is reduced, or as the {\sl highest
root} if the root system is nonreduced.

The choice of basis determines a subset $\RR^+$ of $\RR$ (positive roots);  with the notation
$\RR^-:=\RR^+$ we have $\RR=\RR^+\cup \RR^-$. As usual, $\Q
=\oplus_{i=1}^n\Z\a_i$ denotes the root lattice of $\RR$. Let
$\{\l_i\}_{1\leq i\leq n}$ and $\{\l_i^\vee\}_{1\leq i\leq n}$ be
the fundamental weights, respectively the fundamental coweights
associated to $\RR^+$, and denote by $P=\oplus_{i=1}^n\Z\l_i$ the
weight lattice. An element of $P$ will be called dominant if it is
a linear combination of the fundamental weights with non--negative
integer coefficients. Similarly an anti--dominant weight is a
linear combination of the fundamental weights with non--positive
integer coefficients.

 If $\RR$ is a reduced root system, let $\G$  be a simple complex Lie algebra such that $\Gc$ is a Cartan subalgebra and the associated root system is $\RR$. Also, let $\bo\supset \Gc$ be the Borel subalgebra  determined by $\RR^+$ and let $\bo^-$ be the opposite Borel subalgebra. The simply connected complex algebraic group with Lie algebra $\G$ is denoted by $G$ and $T$, $B$, and $B^-$ denote the subgroups corresponding to $\Gc,\bo$ and $\bo^-$, respectively.

The real vector space $\Gc^*$ has a canonical scalar product
$(\cdot, \cdot)$ which we normalize such that it gives square
length 2 to the short roots in $\RR$ (if there is only one root
length we consider all roots to be short); if $\RR$ is not reduced
 we normalize the scalar product such that the roots have square
length 1, 2 or 4. We will use $\RR_s$ and $\RR_\ell$ to refer to
the short and respectively long roots in $\RR$; if the root system
is nonreduced we will also use $\RR_m$ to refer to the roots of
length 2. We will identify the vector space $\Gc$ with its dual
using this scalar product. Under this identification
$\a^\vee=2\a/(\a,\a)$ for any root $\a$.

To any finite root system as above we associate an {\sl
affine root system} $R$ as follows. Let ${\rm Aff}(\Gc)$ be the space of
affine linear transformations of $\Gc$. As a vector space, it can
be identified to $\Gc^*\oplus \Re\d$ via
$$
(f+c\d)(x)=f(x)+c, \quad \text{for~} f\in \Gc^*, ~x\in\Gc
\text{~and~} c\in \Re
$$
Assume first that $\RR$ is reduced, and let $r$ denote the maximal
number of laces connecting two vertices in its Dynkin diagram.
Define
$$R:=(\RR_s+\Z\d)\cup(\RR_\ell +r\Z\d)\subset \Gc^*\oplus \Re\d$$
If the finite root system $\RR$ is nonreduced define
$$R:=(\RR_s+\frac{1}{2}\Z\d)\cup(\RR_m+\Z\d)\cup(\RR_\ell +\Z\d)$$
Note that in the latter case $R$ is itself a nonreduced root
system. Let us also consider the reduced root systems
$$
\RR_{nd}:=\{\a\in \RR~|~\a/2\not \in \RR\}\quad\text{and}\quad
R_{nd}:=\{\a\in R~|~\a/2\not \in R\}
$$
$$
\RR_{nm}:=\{\a\in \RR~|~2\a\not \in \RR\}\quad\text{and}\quad
R_{nm}:=\{\a\in R~|~2\a\not \in R\}
$$

The set of affine positive roots $R^+$ consists of affine roots of
the form $\a+k\d$ such that $k$ is non--negative if $\a$ is a
positive root, and $k$ is strictly positive if $\a$ is a negative
root. The affine simple roots are $\{\a_i\}_{0\leq i\leq n}$, where
we set $\a_0:=\d-\th$ if $\RR$ is reduced and
$\a_0:=\frac{1}{2}(\d-\th)$ otherwise. In fact, to make some
formulas uniform we set $\a_0:=c_0^{-1}(\d-\th)$, where $c_0$
equals $1$ or 2 depending on whether $\RR$ is reduced or not. If $\a_i$ is a simple root,
then $\a_i^*$ denotes its unique scalar multiple which belongs to $R^+_{nm}$.
Note that $\RR_{nd}$ and $R_{nd}$ have the same basis as $\RR$ and
$R$, respectively. Also, a basis for $\RR_{nm}$ and $R_{nm}$ is given by $\{\a_i^*\}_{1\le i\le n}$ and
$\{\a_i^*\}_{0\le i\le n}$, respectively. The root lattice of $R$ is defined as $Q=\oplus_{i=0}^n \Z\a_i$.

Abstractly, an affine root system is a subset $\Phi\subset {\rm
Aff}(V)$ of the space of affine--linear functions on a real vector
vector space $V$, consisting of non--constant functions which
satisfy the usual axioms for root systems. As in the case of
finite root systems, a classification of the irreducible affine
root systems is available (see, for example, \cite[Section
1.3]{macbook}). The affine root systems $R$ defined above
are just a subset of all the irreducible affine root systems.
However, the configuration of vanishing hyperplanes of elements of
an irreducible affine root system $\Phi$ coincides with the
corresponding configuration of hyperplanes  associated to a {\sl
unique} affine root system $R$ as above. The objects we are concerned with 
here depend in a larger amount on the Weyl group associated
to an affine root system rather that on the root system itself,
and our restriction reflects that. Moreover, the nonreduced affine
root system considered above contains as subsystems all the other
nonreduced irreducible affine root systems and all the reduced
irreducible affine root systems of classical type of the same rank.

%%%%%%%%%%%%%%%%%%%%%%%%%%%%%%%%%%%%%%%%%%%%%%%%%%%%%%%%%%%%%%%%%%%%%%%%%%%%%
%%%%%%%%%%%%%%%%%%%%%%%%%%%%%%%%%%%%%%%%%%%%%%%%%%%%%%%%%%%%%%%%%%%%%%%%%%%%%

\subsection{Affine Weyl groups}\label{awg}

 The scalar product on $\Gc^*$ can be extended to a
non--degenerate bilinear form on the real vector space
$$\G^*:=\Gc^*\oplus \Re\d\oplus \Re\L_0$$ by requiring that $(\d,\Gc^*\oplus \Re\d)=
(\L_0,\Gc^*\oplus \Re\L_0)=0$ and $(\d,\L_0)=1$. Given $\a\in R$
and $x\in \G^*$ let
$$
s_\a(x):=x-\frac{2(x,\a)}{(\a,\a)}\a\
$$
The {\sl affine Weyl group} $W$ is the subgroup of ${\rm
GL}(\G^*)$ generated by all $s_\a$ (the simple reflections
$s_i=s_{\a_i}$ are enough). The {\sl finite Weyl group} $\W$ is
the subgroup generated by $s_1,\dots,s_n$. The bilinear form on
$\G^*$ is equivariant with respect to the affine Weyl group
action. Both the finite and the affine Weyl group are Coxeter
groups and they can be abstractly defined as generated by
$s_1,\dots,s_n$, respectively $s_0,\dots,s_n$, and the following relations:
\begin{enumerate}
\item[a)] reflection relations: $s_i^2=1$; \item[b)] braid
relations: $s_is_j\cdots =s_js_i\cdots $ (there are $m_{ij}$
factors on each side, $m_{ij}$ being equal to $2,3,4,6$ if the
number of laces connecting the corresponding nodes in the Dynkin
diagram is $0,1,2,3$ respectively).
\end{enumerate}

The affine Weyl group could also be presented as a semidirect
product in the following way: it is the semidirect product of $\W$
and the lattice $\Q$ (regarded as an abelian group with elements
$\tau_\mu$, where $\mu$ is in $\Q$), the finite Weyl group acting
on the root lattice as follows
$$
\w\tau_\mu \w^{-1}=\tau_{\w(\mu)}
$$
Since the finite Weyl group  acts on the weight lattice, we
can also consider the {\sl extended Weyl group} $W^e$ defined as
the semidirect product between $\W$ and $P$. Unlike the the affine
Weyl group, $W^e$ is not a Coxeter group. However, $W$ is a normal
subgroup of $W^e$ and the quotient is finite.

For $s$ a real number, $\G^*_s=\{ x\in\G\ ;\ (x,\d)=s\}$ is the
level $s$ of $\G^*$. We have
$$
\G^*_s=\G^*_0+s\L_0=\Gc^*+{\Bbb R}\d+s\L_0\ .
$$
 The action of $W$
preserves each $\G^*_s$ and  we can identify each level
canonically with $\G^*_0$ and obtain an (affine) action
of $W$ on $\G^*_0$. If $s_i\in W$ is a simple reflection, write
$s_i(\cdot)$ for the regular action of $s_i$ on $\G^*_0$ and
$s_i\<\cdot\>$ for the affine action of $s_i$ on $\G^*_0$
corresponding to the level one action. For example, the level zero
action of $s_0$ and $\tau_\mu$ is
\begin{eqnarray*}
s_0(x)     & = & s_\th(x)+(x,\th)c_0^{-1}\d\ ,\\
\tau_\mu(x)   & = & x - (x,\mu)\d\ ,
\end{eqnarray*}
and the level one action of the same elements is
\begin{eqnarray*}
s_0\<x\>   & = & s_\th(x)+(x,\th)c_0^{-1}\d-\a_0\ ,\\
\l_\mu\<x\> & = & x + \mu - (x,\mu)\d -\frac{1}{2}|\mu|^2\d \ .
\end{eqnarray*}
The level one action on $\G^*_0$ induces an affine action of $W$ on $\Gc^*$. As a matter of notation, we write $w\cdot x$ 
for this affine action of $w\in W$ on $x\in \Gc^*$. For example,
\begin{eqnarray*}
s_0\cdot x   & = & s_\th(x)+c_0^{-1}\th\ ,\\
\tau_\mu\cdot x & = & x  + \mu  \ .
\end{eqnarray*}

The fundamental alcove is defined as
\begin{equation}\label{affinechamber}
\cal C:=\{ x\in \Gc^*\ | \ (x+\L_0,\a_i^\vee)\geq 0\ ,\ 0\leq
i\leq n\}
\end{equation}
The non-zero elements of $\cal O_P:=P\bigcap \cal C$ are the
so-called minuscule weights. Note that each orbit of the
affine action of $W$ on $P$ contains either the origin or a unique
fundamental weight $\l_i$ (to keep the notation
consistent we set $\l_0=0$).
%In fact, for any intermediary lattice
%$\Q\subset \L\subset P$ on which the affine Weyl group acts, the
%elements of $\cal O_\L:=\L\cap \cal C$ form a set of
%representatives for the orbits of $W$ on $\L$.
If we denote $$\Omega:=\{w\in W^e~|~ w\cdot \cal C=\cal C\} $$
then, $W^e=\Omega\ltimes W$. The group $\Omega$ is finite, of
order the size of $\cal O_P$. In fact, we can parameterize
$\Omega$ by the elements of $\cal O_P$ as follows: for each $\l\in
\cal O_P$ let $\omega_\l$ denote the unique element of $\Omega$
for which $\omega_\l(0)=\l$. It is easy to check that
$\omega_\l=\tau_\l\w_\l$. For the definition of $\w_\l$ see Section \ref{coset representatives} below.

If we examine the orbits of the level zero action of the affine Weyl group $W$  on the affine root system $R$ we find the following: 
\begin{enumerate}
\item[a)] if
$\RR$ is reduced there are precisely as many orbits as root
lengths;
\item[b)] if $\RR$ is nonreduced of rank at least two, then there are five orbits:
$$
W(2\a_0)=\RR_\ell+2\Z\d+\d, \quad
W(\a_0)=\RR_s+\Z\d+\frac{1}{2}\d,\quad W(\a_1)=\RR_m+\Z\d
$$
$$
 W(2\a_n)=\RR_\ell+2\Z\d\quad\text{and}\quad W(\a_n)=\RR_s+\Z\d
$$
\item[c)] if $\RR$ is nonreduced of rank one then there are only for orbits: $W(2\a_0)$, 
$W(\a_0)$, $W(2\a_1)$ and $W(\a_1)$. 
\end{enumerate}
%%%%%%%%%%%%%%%%%%%%%%%%%%%%%%%%%%%%%%%%%%%%%%%%%%%%%%%%%%%%%%%%%%%%%%
%%%%%%%%%%%%%%%%%%%%%%%%%%%%%%%%%%%%%%%%%%%%%%%%%%%%%%%%%%%%%%%%%%%%%%
%%%%%%%%%%%%%%%%%%%%%%%%%%%%%%%%%%%%%%%%%%%%%%%%%%%%%%%%%%%%%%%%%%%%%%

\subsection{The length function}
For each $w$ in $W$ let $\ell(w)$ be the length of a reduced (i.e.
shortest) decomposition of $w$ in terms of simple reflections. The
length of $w$ can be also geometrically described as follows. For
any affine root $\a$, denote by $H_\a$ the affine hyperplane
consisting of fixed points of the affine action of $s_\a$ on
$\Gc^*$. Then, $\ell(w)$ equals the number of affine hyperplanes
$H_\a$ separating $\cal C$ and $w\cdot\cal C$. Since the affine
action of $W^e$ on $\Gc^*$ preserves the set 
$\{H_\a\}_{\a\in R}$, we can use the geometric point of view to
define the length of any element of $W^e$. For example, the
elements of $\Omega$ will have length zero.

For $w$ in $W$ we have $\ell(w)=|\Pi(w)| $, where $
\Pi(w)=\{\a\in R_{nd}^+ ~ |~w(\a)\in R_{nd}^-\}$. We denote
 $ {}^c\Pi(w):=\{\a\in R_{nd}^+\ |\ w(\a)\in R_{nd}^+\}$. If
$w=s_{j_p}\cdots s_{j_1}$ is a reduced decomposition, then
\begin{equation}\label{Pi}
\Pi(w)=\{\a^{(i)}\ |\ 1\leq i\leq p\},
\end{equation}
 with
$\a^{(i)}=s_{j_1}\cdots s_{j_{i-1}}(\a_{j_i})$. An easy
check shows that for any $w$ in $W$ we have
\begin{equation}\label{Piw's}
w^{-1}\left(\Pi(w^{-1})\right)=-\Pi(w) \ \text{ and }\
w^{-1}\left({}^c\Pi(w^{-1})\right)={}^c\Pi(w)
\end{equation}
The following formula is well-known (see, for example,
\cite{lusztig}). If $\w\in\W$ and $\l\in P$, then
\begin{equation}\label{lusztiglength}
\ell(\w\tau_\l)=\sum_{\a\in\Pi(\w)} |(\l,\a^\vee) +1|
+\sum_{\a\in{}^c\Pi(\w) } |(\l,\a^\vee) |
\end{equation}
Let us derive a few immediate consequences which will be useful
later on.
\begin{Lm}\label{length}
Assume that $\l$ and $\mu$ are dominant weights  and that $\w$ is
an element of $\W$. Then, the following equalities hold:
\begin{enumerate}
\item $\ell(\tau_{\l+\mu})=\ell(\tau_\l)+\ell(\tau_\mu) $ \item
$\ell(\w\tau_\l)=\ell(\w)+\ell(\tau_\l)$ \item
$\ell(\tau_{\w(\l)})=\ell(\tau_\l)$
\end{enumerate}
\end{Lm}
\begin{proof}
 The first two claims follow directly from the formula
 (\ref{lusztiglength}) if we keep in mind that the scalar product $(\l,\a^\vee)$
 is non-negative if $\l$ is dominant and $\a$ is a positive finite
 root. To prove the third statement note that
 \begin{eqnarray*}
\ell(\tau_\l) &=& \sum_{\a\in \RR^+_{nd}} |(\l,\a^\vee)|\\
&=& \sum_{\a\in \Pi(\w)} |(\l,\a^\vee)| +\sum_{\a\in {}^c\Pi(\w)}
|(\l,\a^\vee)|\\
&=& \sum_{\a\in \Pi(\w^{-1})} |(\l,\w^{-1}(\a^\vee))| +\sum_{\a\in
{}^c\Pi(\w^{-1})} |(\l,\w^{-1}(\a^\vee))|\\
&=& \sum_{\a\in \Pi(\w^{-1})} |(\w(\l),\a^\vee)| +\sum_{\a\in
{}^c\Pi(\w^{-1})} |(\w(\l),\a^\vee)\\
&=& \ell(\tau_{\w(\l)})
 \end{eqnarray*}
 Along the way we have used the equalities (\ref{Piw's}).
\end{proof}

%%%%%%%%%%%%%%%%%%%%%%%%%%%%%%%%%%%%%%%%%%%
%%%%%%%%%%%%%%%%%%%%%%%%%%%%%%%%%%%%%%%%%%%
%%%%%%%%%%%%%%%%%%%%%%%%%%%%%%%%%%%%%%%%%%%
\subsection{Coset representatives}\label{coset representatives}
For each weight $\l$ define $\l_-$, respectively $\tilde\l$, to be
the unique element in $\W(\l)$, respectively $W\cdot\l$, which is
an anti-dominant weight, respectively an element of $\cal O_P$
(i.e. either zero or a minuscule weight), and
$\w_\l\hspace{-0.2cm}^{-1}\in\ \W$, respectively $w_\l^{-1}\in W$,
to be the unique minimal length elements by which this is achieved.
Also, for each weight $\l$ define $\l_+$ to be the unique dominant element
in $\W(\l)$ and denote by $w_\circ$ the maximal
length element in $\W$. 

Clearly, the element $\w_\l$ is the minimal length representative in its  left coset $\w_\l {\rm Stab}_{\W}(\l_-)\subset \W$. The element $w_\l$ can be equivalently
described as the minimal length representative of the coset
$\tau_\l\W\omega_{\tilde \l}^{-1}\subset W$. Similarly, we
consider $v_\l$, the unique maximal length representative of the
coset $\tau_\l\W\omega_{\tilde \l}^{-1}=w_\l\omega_{\tilde
\l}\W\omega_{\tilde \l}^{-1}$. In fact, the group $\omega_{\tilde
\l}\W\omega_{\tilde \l}^{-1}$ is the stabilizer  ${\rm Stab}_W(\tilde \l)$ which will be 
denoted by $W_{\tilde \l}$. Its maximal
length element is $w_{\circ,\tilde\l}:=\omega_{\tilde
\l}w_\circ\omega_{\tilde \l}^{-1}$ and $v_\l$ and $w_\l$ are
related by the formula
\begin{equation}\label{eq12}
v_\l=w_\l w_{\circ,\tilde\l}
\end{equation}
Moreover,
\begin{equation}\label{eq14}
\ell(v_\l)=\ell(w_\l)+\ell(w_{\circ,\tilde\l})=\ell(w_\l)+\ell(w_{\circ})
\end{equation}

Let us recall from \cite{ion} the following result.
\begin{Lm}
With the above notation
\begin{enumerate}\label{lemma2}
\item $ \Pi(\w_\l\hspace{-0.2cm}^{-1})=\{\a\in \RR_{nd}^+\ |\
(\l,\a)>0 \} $ \item $ \Pi(w_\l^{-1})=\{\a\in R_{nd}^+\ |\
(\l+\L_0,\a)<0 \} $
\end{enumerate}
\end{Lm}
The following technical result will be used later.
\begin{Lm}\label{lemma5}
Let $\l$ be a weight and let $\b$ be a root in $\RR$ such
that $\a=\b+k\d$ is a positive affine root. If $(\a,\l+\L_0)<0$
then $\w_\l^{-1}(\b)$ belongs to $\RR^+$.
\end{Lm}
\begin{proof} Let us remark that it is enough to prove our
result for some positive scaling of the root $\a$ and therefore we
can safely assume that $\a\in R^+_{nd}$.

Since $\a=\b+k\d$ is a positive affine root we have to analyze two
possible situations. First, assume that $\b\in\RR^+$ and $k\geq
0$. In this case, $(\a,\l+\L_0)<0$ implies that $(\b,\l)<0$ and
the above Lemma tells us that $\w_\l^{-1}(\b)\in \RR^+$.

The other possible situation is $\b\in\RR^-$ and $k> 0$. In this
case, $(\a,\l+\L_0)<0$ implies that $(-\b,\l)>0$ and since $-\b\in
\RR^+$ we obtain that $\w_\l^{-1}(-\b)\in \RR^-$ or, equivalently,
$\w_\l^{-1}(\b)\in \RR^+$.
\end{proof}

%%%%%%%%%%%%%%%%%%%%%%%%%%%%%%%%%%%%%%%%%%%%%%%%%%%%%%%%%%%%%%%%%%%%%
%%%%%%%%%%%%%%%%%%%%%%%%%%%%%%%%%%%%%%%%%%%%%%%%%%%%%%%%%%%%%%%%%%%%%
%%%%%%%%%%%%%%%%%%%%%%%%%%%%%%%%%%%%%%%%%%%%%%%%%%%%%%%%%%%%%%%%%%%%%

\subsection{The Bruhat order}

The Bruhat order is a partial order on any Coxeter group defined
in way compatible with the length function. For an element $w$ and a root $\a$ we
write $w<s_\a w$\  if and only if \ $\ell(s_\a w)=\ell(w)+1$. The
transitive closure of this relation is called the Bruhat order on $W$.
The terminology is motivated by the way this order arises for
Weyl groups in connection with inclusions among the closures of the Bruhat cells
 of a corresponding connected, simple algebraic group. For the
basic properties of the Bruhat order we refer to Chapter 5 in
\cite{humphreys}. Let us list a few of them (the first two
properties completely characterize the Bruhat order):
\begin{enumerate}
\item For each $\a\in R^+$ we have $s_\a w<w$ if and only if $\a$
is in $\Pi(w^{-1})$ ; \item $w'< w$ if and only if $w'$ can be
obtained by omitting some factors in a reduced decomposition
of $w$ ; \item If $s_i$ is a simple reflection and $w' \leq w$
then either $s_i w' \leq w$ or $s_i w' \leq s_iw$ (or both). For
example, if $\ell(s_iw')-\ell(w')\leq\ell(s_iw)-\ell(w)$ then $s_i w'
\leq s_iw$.
\end{enumerate}
We can use the Bruhat order on $W$ to define a partial order on
the weight lattice which will also be called the Bruhat order. For
any $\l,\mu\in P$ we write
\begin{equation}\label{order}
\l<\mu\ \ \text{if and only if}\ \tilde \l=\tilde \mu\  \text{and}\ w_\l<w_\mu
\end{equation}
The minimal elements of $P$ with respect to this
partial order are the minuscule weights. The next result shows that
this partial order relation could have been defined as well using
the elements $v_\l$ instead of $w_\l$.
\begin{Lm}
Let $\l$ and $\mu$ be two  weights. Then $w_\mu<w_\l$ if and only
if $v_\mu<v_\l$.
\end{Lm}
\begin{proof}
Straightforward from  (\ref{eq12}) and the third property of the
Bruhat order.
\end{proof}
\begin{Lm}\label{lemma6}
Let $\l$ be a weight. We have
$$
\{x\in W~|~x\leq v_\l\}=\bigcup_{\mu\leq\l} v_\mu W_{\tilde\l}
$$
\end{Lm}
\begin{proof}
If $\mu$ is a weight such that $\mu\leq\l$ then, by the above
Lemma, $v_\mu\leq v_\l$. Since $v_\mu$ is the maximal element of
the coset $v_\mu W_{\tilde\l}$, we obtain that $y\leq v_\l$ for
any element $y$ in $v_\mu W_{\tilde\l}$.

Conversely, let $x\in W$ such that $x\leq v_\l$. The third
property of the Bruhat order together with the definition of
$v_\l$ imply that $z\leq v_\l$ for any $z\in x W_{\tilde\l}$. The
left coset $x W_{\tilde\l}$ is of the form $v_\mu W_{\tilde\l}$
for some weight $\mu$ for which $\tilde\mu=\tilde\l$. Therefore, as claimed,
we obtain that $x\in v_\mu W_{\tilde\l}$ and $v_\mu\leq v_\l$.
\end{proof}

The following result can be found in \cite{ion}.
\begin{Lm}\label{lemma3}
Let $\l$ be a weight and $\a_i$ be a simple affine root such that
$s_i\cdot\l\not =\l$. The following statements hold:
\begin{enumerate}
\item We  have, $\w_{s_i\cdot\l}=s_i\w_\l$, unless $i=0$, in which
case $\w_{s_0\cdot\l}=s_\th \w_\l$. 
\item Moreover,  $s_i\cdot\l>\l$ if and only if
$(\a_i,\l+\L_0)> 0$. In particular $\l_-$, respectively $\l_+$,
are the maximal element, respectively the minimal element in
$\W(\l)$ with respect to the Bruhat order.
\end{enumerate}
\end{Lm}
We close this section with a consequence of the above result.
\begin{Lm}\label{lemma4}
 Let $\l$ be an anti--dominant  weight and $\mu\in \W(\l)$. Then
 \begin{enumerate}
 \item $w_{\l}=\w_\mu^{-1}w_\mu$ and
$\ell(w_{\l})=\ell(w_\mu)+\ell(\w_\mu)$ 
\item $w_\l\omega_{\tilde \l}=\tau_\l$
\item $\tau_\mu\w_\mu=w_\mu\omega_{\tilde \mu}$ and
$\ell(\tau_\mu)=\ell(w_\mu)+\ell(\w_\mu)$
\end{enumerate}
\end{Lm}
\begin{proof} (1)
Let us note  that if we fix a reduced decomposition $s_{j_p}\cdots
s_{j_1}$ for $\w_{\mu}^{-1}$ then
\begin{equation}\label{eq1}
\l=s_{j_p}\cdots s_{j_1}(\mu) > \cdots > s_{j_2}s_{j_1}(\mu)
> s_{j_1}(\mu)
>\mu
\end{equation}
Indeed, for any $1\leq i\leq p$ we have
$$(\a_{j_i},s_{j_{i-1}}\cdots s_{j_1}(\mu))=(s_{j_{1}}\cdots
s_{j_{i-1}}(\a_{j_i}),\mu)$$ and from  equation (\ref{Pi}) we
know that $s_{j_{1}}\cdots s_{j_{i-1}}(\a_{j_i})$ belongs to
$\Pi(\w_\mu^{-1})$. Furthermore, Lemma \ref{lemma2} implies that
$(s_{j_{1}}\cdots s_{j_{i-1}}(\a_{j_i}),\mu)>0$ and Lemma
\ref{lemma3} immediately gives us (\ref{eq1}). We conclude that
$w_{\l}=\w_\mu^{-1}w_\mu$ and
$\ell(w_{\l})=\ell(w_\mu)+\ell(\w_\mu)$.

(2) By definition, $w_{\l}$ is the unique minimal length element
in the coset $\tau_{\l}\W\omega_{\tilde \l}^{-1}$. Since
$\omega_{\tilde \l}$ has length zero it is enough to show that
$\ell(\tau_{\l}\w)\geq \ell(\tau_{\l})$ for all $\w\in \W$.
Indeed, $\ell(\tau_{\l}\w)=\ell(\w^{-1}\tau^{-1}_{\l})$ and since
the element $-\l$ is dominant ($\l$ \ being anti--dominant) Lemma
\ref{length} implies the desired result.

(3) The statement follows immediately from (1), (2) and Lemma
\ref{length}.
\end{proof}

%%%%%%%%%%%%%%%%%%%%%%%%%%%%%%%%%%%%%%%%%%%%%%%%%%%%%%%%%%%%%%%%%%%%%%%%%%%%
%%%%%%%%%%%%%%%%%%%%%%%%%%%%%%%%%%%%%%%%%%%%%%%%%%%%%%%%%%%%%%%%%%%%%%%%%%%%
%%%%%%%%%%%%%%%%%%%%%%%%%%%%%%%%%%%%%%%%%%%%%%%%%%%%%%%%%%%%%%%%%%%%%%%%%%%%

\section{Nonsymmetric Macdonald polynomials}\label{section2}

\subsection{Parameters and conventions}
Let us introduce a field $\F$ (of parameters) as follows. Let
$t=(t_\a)_{\a\in R}$ be a set of parameters which is indexed by
the set of affine roots and has the property that $t_\a=t_\b$ if
and only if the affine roots $\a$ and $\b$ belong to the same
orbit under the action of $W$ on $R$. It will be convenient to have also the following
convention: if $\a$ is not an affine root then $t_\a=1$. Let $q$
be another parameter and let $m$ be the lowest common denominator
of the rational numbers $\{(\a_j,\l_k)\ |\ 1\leq j,k\leq n \}$.
The field $\F=\F_{q,t}$ is defined as the field of rational
functions in $q^{\frac{1}{m}}$ and
$t^{\frac{1}{2}}=(t_\a^{\frac{1}{2}})_{\a\in R}$ with rational coefficients. We will also use
the field of rational functions in
$t^{\frac{1}{2}}=(t_\a^{\frac{1}{2}})_{\a\in R}$ denoted by
$\F_t$.

As it follows from the discussion at the end of Section \ref{awg} there are only a small number of distinct parameters.
If the root system $R$ is reduced then  there are as many distinct
parameters $t_\a$ as root lengths: at most two, which we denote by $t_s$ (the one corresponding to short roots) and $t_\ell$ (the one 
corresponding to long roots). In this
case, to avoid unnecessary notational complexity we use
$t_i$ to refer to the parameter $t_{\a_i}$ corresponding to the affine simple root $\a_i$. 

If $R$ is nonreduced then the action of the affine Weyl group on
the affine root system has five orbits
$W(2\a_0),~W(\a_0),~W(\a_n), ~W(2\a_n)$ and $W(\a_1)$ (note that
the last orbit is empty if $R$ has rank one) and we denote the
corresponding parameters by $t_{01},~ t_{02},~ t_{03},~ t_n$ and
$t:=t_1=\cdots=t_{n-1}$, respectively. The relation with the notation used in \cite{sahi} is the
following: $t_0,~u_0,~u_n$ used in \cite{sahi} are respectively $t_{01},~t_{02},~t_{03}$ in our notation.

To avoid distinguishing among the reduced and nonreduced case later on we find convenient to define 
$t_{01}=t_{02}=t_{03}:=t_0$ in the reduced case.

%%%%%%%%%%%%%%%%%%%%%%%%%%%%%%%%%%%%%%%%%%%%%%%%%%%%%%%%%%%%%%
%%%%%%%%%%%%%%%%%%%%%%%%%%%%%%%%%%%%%%%%%%%%%%%%%%%%%%%%%%%%%%
%%%%%%%%%%%%%%%%%%%%%%%%%%%%%%%%%%%%%%%%%%%%%%%%%%%%%%%%%%%%%%

\subsection{Double affine Hecke algebras}

The algebra $\cal R=\F[e^\l;\l\in P]$ is the group
$\F$-algebra of the lattice $P$. Similarly, the algebra $\cal
R_t=\F_t[e^\l;\l\in P]$ is the group $\F_t$-algebra of the lattice
$P$. In the discussion that follows we refer to the following group
$\F$--algebras of the root lattice: $\cal{Q}_Y:=\F[Y_\mu;\mu\in
\Q]$ and $\cal{Q}_X:=\F[X_\b;\b\in \Q]$. We will also use the
following notation: for $\mu\in \Q$ and $k\in \frac{1}{m}\Z$ let
$e^{\mu+k\d}:=q^{-k}e^\mu$, $X_{\b+k\d}:=q^{-k}X_\b$ and
$Y_{\mu+k\d}:=q^kY_{\mu}$.

In the reduced case, the double affine Hecke algebras were
introduced by Cherednik (see, for example, \cite{c1}) in his work
on affine quantum Knizhnik--Zamolodchikov equations and on
Macdonald's conjectures. In the nonreduced case the definition is
due to Sahi \cite{sahi}. We give here the symmetric definition of
the double affine Hecke algebras obtained in \cite{is}.
\begin{Def}
The double affine Hecke algebra $\tilde \H$ associated to the root
system $\RR$ is the $\F$--algebra described by generators and
relations as follows:

\underline{Generators}: One generator $T_i$ for each simple root
$\a_i$, with the exception of the affine simple root $\a_0$ for
which we associate three generators $T_{01}$, $T_{02}$ and $T_{03}$.

\underline{Relations}: a) Each pair of generators satisfies the
same braid relations as the corresponding pair of simple reflections.

\hspace{1.52cm} b) If  there is a simple root $\a$ such that
$(\a,\a_0^\vee)=-2$ then the following relation also holds
$$
T_{01}T_\a^{-1}T_{03}T_\a=T_\a^{-1}T_{03}T_\a T_{01}
$$

\hspace{1.52cm} c) The quadratic relations
\begin{eqnarray*}
T_i^2 &=& (t_i^{\frac{1}{2}} -t_i^{-\frac{1}{2}})T_i +1, \quad \ \
\text{for all } \ \ 1\leq i\leq n, \text{and} \\
T_{0j}^2 &=& (t_{0j}^{\frac{1}{2}}
-t_{0j}^{-\frac{1}{2}})T_{0j}+1, \ \ \ \text{for} \ \ 1\leq j\leq
3
\end{eqnarray*}

\hspace{1.52cm} c) The relation $$
T_{01}T_{02}T_{03}T_{s_\th}=q^{-c_0^{-1}} $$
\end{Def}
In the case of a reduced root system the quadratic relations for
the elements $T_{0j}$ need not be imposed, since they are a
consequence of the other relations. However it is absolutely
necessary to impose them for nonreduced root systems. For
nonreduced root systems,  the relationship between the generators
$T_{0j}$ and  the notation used in \cite{sahi}  is the
following: $T_0,~U_0,~U_n$ used in \cite{sahi} are respectively $T_{01},~T_{02},~T_{03}$ in our notation.

The elements $T_1, \dots, T_n$ generate the  {\sl finite Hecke
algebra} $\h$. There are countably many copies of the {\sl affine
Hecke algebra} associated to the affine root system $R$ inside
$\tilde \H$; we will distinguish only two of them: $\H_Y$ which is
the subalgebra generated by $T_{01},T_1,\cdots, T_n$, and $\H_X$
which is the subalgebra generated by $T_{03},T_1,\cdots, T_n$.
There are natural bases of $\H_X$, $\H_Y$ and $\h$: $\{T_w\}_w$
indexed by $w$ in $W$ and in $\W$ respectively, where
$T_w=T_{i_l}\cdots T_{i_1}$ if $w=s_{i_l}\cdots s_{i_1}$ is a
reduced  expression of $w$ in terms of simple reflections. Let us
recall the well--known result of Bernstein (unpublished) and
Lusztig \cite{lusztig} on the structure of affine Hecke algebras
as it applies to $\H_Y$ and $\H_X$.

\begin{Prop}\label{HXHY} With the above notation we have
\begin{enumerate}
\item The affine Hecke algebra $\H_Y$ is generated by the finite
Hecke algebra $\h$ and the group algebra $\cal Q_Y$ such that the
following relations are satisfied for any $\mu$ in the root
lattice and any $1\leq i\leq n$\ :
$$\hspace*{-2.8cm}
Y_\mu T_i-T_iY_{s_i(\mu)} = (t_i^{\frac{1}{2}}-t_i^{-\frac{1}{2}})
\frac{Y_\mu-Y_{s_i(\mu)}}{1-Y_{\a_i}}~~ \text{ if }~ 2\a_i\not \in
\RR
$$
$$ Y_\mu T_n  - T_nY_{s_n(\mu)}  =\left(
t_n^\frac{1}{2}-t_n^{-\frac{1}{2}}+( t_{01}^{\frac{1}{2}}-
t_{01}^{-\frac{1}{2}})Y_{\a_n}\right)\frac{Y_\mu-Y_{s_n(\mu)}}{1-Y_{2\a_n}}~
~\text{ if }~~~ 2\a_n\in \RR
$$
In this description $Y_{-c_0^{-1}\th } =T_{s_\th}T_{01}$.

\item The affine Hecke algebra $\H_X$ is generated by the finite
Hecke algebra $\h$ and the group algebra $\cal Q_X$ such that the
following relations are satisfied for any $\mu$ in the root
lattice and any $1\leq i\leq n$\ :
$$\hspace*{-3.2cm} T_iX_\mu-X_{s_i(\mu)}T_i = (t_i^{\frac{1}{2}}-t_i^{-\frac{1}{2}})
\frac{X_\mu-X_{s_i(\mu)}}{1-X_{-\a_i}}~~ \text{ if }~ 2\a_i\not
\in \RR
$$
$$ T_n X_\mu - X_{s_n(\mu)}T_n  =\left(
t_n^\frac{1}{2}-t_n^{-\frac{1}{2}}+( t_{03}^{\frac{1}{2}}-
t_{03}^{-\frac{1}{2}})X_{-\a_n}\right)\frac{X_\mu-X_{s_n(\mu)}}{1-X_{-2\a_n}}~
~\text{ if }~~~ 2\a_n\in \RR
$$
In this description $X_{c_0^{-1}\th } =T_{03}T_{s_\th}$.
\end{enumerate}
\end{Prop}
\begin{Rem}\label{remarkT02} We note that with the above notation
$$T_{02}=T^{-1}_{01}X_{\a_0}=Y_{-\a_0}T_{03}^{-1}=q^{-c_0^{-1}}Y_{c_0^{-1}\th}T_{s_\th}X_{-c_0^{-1}\th} $$
Therefore, $T_{01}=T_0$ and $T_{02}=T_{\<0\>}$ with the notation
used in \cite{ion}.
\end{Rem}

To define an action of $\tilde \H$ one needs only to define the
action of the generators $T_i$, $1\leq i\leq n$ and $T_{0j}$,
$1\leq j\leq 3$. However, from Remark \ref{remarkT02} and
Proposition \ref{HXHY} it is clear that we can equivalently define
a representation of the double affine Hecke algebra by only
specifying the action of $T_{01}$, $T_i$, $1\leq i\leq n$ and
$\cal Q_X$. From the work of Cherednik (in the reduced case) and
Sahi (in the nonreduced case) we know that the following formulas
define a faithful representation of $\tilde \H$ on $\cal R$
\begin{eqnarray*}
X_\mu\cdot e^\l &=& e^{\l+\mu}\ \ \text{for} \ \mu\in \ \Q\\
 T_i\cdot
e^\l &=& t_i^{\frac{1}{2}}e^{s_i(\l)}
+(t_i^{\frac{1}{2}}-t_i^{-\frac{1}{2}})
\frac{e^\l-e^{s_i(\l)}}{1-e^{-\a_i}}\ \ \text{if}\ 1\leq i\leq n
~\text{and} ~2\a_i\not\in \RR\\
T_n\cdot e^\l&=& t_n^{\frac{1}{2}}e^{s_n(\l)}+\left(
t_n^\frac{1}{2}-t_n^{-\frac{1}{2}}+( t_{03}^{\frac{1}{2}}-
t_{03}^{-\frac{1}{2}})e^{-\a_n}\right)\frac{e^\l-e^{s_n(\l)}}{1-e^{-2\a_n}}~
~\text{ if }~~~ 2\a_n\in \RR\\
 T_{01}\cdot
e^\l &=& t_0^{\frac{1}{2}}e^{s_0(\l)}
+(t_0^{\frac{1}{2}}-t_0^{-\frac{1}{2}})
\frac{e^\l-e^{s_0(\l)}}{1-e^{-\a_0}}\ \ \text{if}\ ~ ~2\a_0\not\in R\\
T_{01}\cdot e^\l&=& t_{01}^{\frac{1}{2}}e^{s_0(\l)}+\left(
t_{01}^\frac{1}{2}-t_{01}^{-\frac{1}{2}}+( t_{02}^{\frac{1}{2}}-
t_{02}^{-\frac{1}{2}})e^{-\a_0}\right)\frac{e^\l-e^{s_0(\l)}}{1-e^{-2\a_0}}~
~\text{ if }~~~ 2\a_0\in R
\end{eqnarray*}
Using Remark \ref{remarkT02} the action of $T_{02}$ is easily
computable. Let us list the results
\begin{eqnarray*} T_{02}\cdot e^\l &=&  t_{0}^{\frac{1}{2}}e^{s_0\<\l\>} +
(t_{0}^{\frac{1}{2}}-t_{0}^{-\frac{1}{2}})
\frac{e^\l-e^{s_0\<\l\>}}{1-e^{-\a_0}}\ \ \text{if}\ ~ ~2\a_0\not\in R\\
T_{02}\cdot e^\l &=& t_{01}^{\frac{1}{2}}e^{s_0\<\l\>} +\left(
t_{02}^\frac{1}{2}-t_{02}^{\frac{1}{2}}+( t_{01}^{\frac{1}{2}}-
t_{01}^{-\frac{1}{2}})e^{-\a_0}\right)\frac{e^\l}{1-e^{-2\a_0}}\\
& &-\left( t_{01}^\frac{1}{2}-t_{01}^{\frac{1}{2}}+(
t_{02}^{\frac{1}{2}}-
t_{02}^{-\frac{1}{2}})e^{-\a_0}\right)\frac{e^{s_0\<\l\>}}{1-e^{-2\a_0}}\
\ \text{if}\ ~ ~2\a_0\in R
\end{eqnarray*}

We also need to consider the extended affine Hecke algebra
$\H_X^e$ which is defined as the semidirect product of $\Omega$
and $\H_X$. The action of $\Omega$ on $\H_X$ is induced from the
action of $\Omega$ on the affine Weyl group: if $\omega\in \Omega$
and $\w\in W$ then $\omega T_w \omega^{-1}=T_{\omega
w\omega^{-1}}$. If we use the notation $T_{w\omega}:= T_w \omega$,
a basis for $\H_X^e$ is given by $\{T_w\}_{w\in W^e}$. The action
of $\H_X$ on $\cal R$ described above can be extended to an action
of $\H_X^e$ by defining
$$
\omega_{\l}\cdot e^\mu=e^\l T_{\w_\l^{-1}}^{-1}\cdot e^{\mu}
$$
for any $\l$ in $\cal O_P$. It is important to note that for any
dominant weights $\nu_1$ and $\nu_2$ the element
$X_{\nu_1-\nu_2}:=T_{\tau_{\nu_1}}T_{\tau_{\nu_2}}^{-1}$ acts on $\cal
R$ as multiplication by $e^{\nu_1-\nu_2}$.

%%%%%%%%%%%%%%%%%%%%%%%%%%%%%%%%%%%%%%%%%%%%%%%%%%%%%%%%%%%%%%%%%%%%%%%%%%%%%
%%%%%%%%%%%%%%%%%%%%%%%%%%%%%%%%%%%%%%%%%%%%%%%%%%%%%%%%%%%%%%%%%%%%%%%%%%%%%
%%%%%%%%%%%%%%%%%%%%%%%%%%%%%%%%%%%%%%%%%%%%%%%%%%%%%%%%%%%%%%%%%%%%%%%%%%%%%

\subsection{The Cherednik scalar product}
The involution of $\F$ which inverts each of the parameters $q,~\{t_\a\}_{\a\in R}$
extends to an involution $\overline{~\cdot~}$ on the algebra $\cal
R$ which sends each $e^\l$ to $e^{-\l}$. Following Cherednik let
us define
\begin{equation}
K(q,t)=\prod_{\a\in R_{nm}^+} \frac{1-e^{\a}}
{(1-t_\a^{-\frac{1}{2}}t_{\frac{\a}{2}}^{-\frac{1}{2}}
e^{\frac{\a}{2}})
(1+t_\a^{-\frac{1}{2}}t_{\frac{\a}{2}}^{\frac{1}{2}}
e^{\frac{\a}{2}})}
\end{equation}
which should be seen as a formal series in the elements $e^\a$
with coefficients in $\Z[t^{-1}][[q^{-1}]]$ (power series in
$q^{-1}$ with coefficients polynomials in $t^{-1}$). Recall that
we agreed to set $t_{\frac{\a}{2}}=1$ if $\a/2$ is not a root. If
$K_0$ denotes the coefficient of $e^0$, then the Cherednik kernel
is by definition
$$
C(q,t):=\frac{K(q,t)}{K_0}
$$
It is a formal series in the elements $e^\a$ with coefficients
rational functions in $q$ and $t$ (see, for example, {\rm
\cite[(5.1.10)]{macbook}}). Moreover, $C(q,t)$ it is fixed by the
above involution.

A scalar product on $\cal R$ can be defined as follows
$$
\<f,g\>_{q,t}:=CT(f\bar g C(q,t))
$$
where $CT(\cdot)$ denotes the constant term (i.e. the
coefficient of $e^0$) of the expression inside the parenthesis.
The scalar product is Hermitian with respect to the involution
$\overline{~\cdot~}$~:
$$\<g,f\>_{q,t}=\overline{\<f,g\>}_{q,t}$$
and the above representation becomes unitary with respect to
$\<\cdot,\cdot\>_{q,t}$.
%%%%%%%%%%%%%%%%%%%%%%%%%%%%%%%%%%%%%%%%%%%%%%%%%%%%%%%%%%%%%%%%%%%%%%%%%%%%%
%%%%%%%%%%%%%%%%%%%%%%%%%%%%%%%%%%%%%%%%%%%%%%%%%%%%%%%%%%%%%%%%%%%%%%%%%%%%%
%%%%%%%%%%%%%%%%%%%%%%%%%%%%%%%%%%%%%%%%%%%%%%%%%%%%%%%%%%%%%%%%%%%%%%%%%%%%%
\subsection{Macdonald polynomials} If $\gamma$ is an element of $\Gc^*\oplus \frac{\d}{c_0}\Z$
we denote by ${\bf q}^{(\gamma,\overline\l)}$  the element of $\F$
$$
q^{(\gamma,\l+\L_0)}(t_{n}^*t_n)^{-\frac{1}{2}(\gamma,
\w_\l(\l_n^\vee))} \prod_{i=1}^{n-1}t_i^{-(\gamma,
\w_\l(\l_i^\vee))}
$$
where $t_n^*$ equals $t_n$ if $R$ is reduced or $t_{01}$ if $R$ is
nonreduced. Under the same conditions let $${\bf
t}^{(\gamma,\overline\l)}:=(t_{n}^*t_n)^{\frac{1}{2}(\gamma,
\w_\l(\l_n^\vee))} \prod_{i=1}^{n-1}t_i^{(\gamma,
\w_\l(\l_i^\vee))}$$ In particular, we have
$q^{(\gamma,\l+\L_0)}={\bf q}^{(\gamma,\overline\l)}{\bf
t}^{(\gamma,\overline\l)}$.

\noi For each $\l\in P$ we can construct a $\F$-algebra morphism
${\rm ev}(\l):\cal Q_Y\to \F$, which sends $Y_\mu$ to ${\bf
q}^{(\mu,\overline \l)}$. If $f(Y)$ is an element of $\cal Q_Y$ we
will write $f(\l)$ for ${\rm ev}(\l)(f)$.

For every weight $\l$ define
$$
\cal R_{\l}=\{f\in\cal R\ | \ Y_\mu \cdot f={\bf
q}^{(\mu,\overline \l)}f \ \text{for any } \mu\in\Q \}.
$$
\begin{Def}
Given a weight $\l$  the nonsymmetric Macdonald polynomial
$E_\l(q,t)$ is the unique element in $\cal R_\l$ in which the
coefficient of $e^\l$ is $1$. If $k\in \Z$ denote
$$E_{\l+k\d}(q,t):=q^{-k}E_\l(q,t)$$
\end{Def}

The nonsymmetric Macdonald polynomials form a basis of $\cal R$
orthogonal with respect to the scalar product $\<\cdot,
\cdot\>_{q,t}$. They are also triangular with respect to the
Bruhat order on the weight lattice. Since the minimal elements for
this order relation are the minuscule weights we immediately
obtain the following
\begin{Prop}\label{Eminuscule}
If $\l$ is a minuscule weight, then $E_\l(q,t)=e^\l$.
\end{Prop}

For any anti--dominant weight $\l$ we write $\cal R^\l$ for the
subspace of $\cal R$ spanned by  $\{E_\mu\ |\ \mu\in\
\W(\l) \}$. The relationship with the symmetric
Macdonald polynomials is the  following.
\begin{Def}
Given an anti--dominant weight $\l$ the symmetric Macdonald
polynomial $P_\l(q,t)$ can be characterized as the unique $\W$-
invariant element in $\cal R^\l$ for which the coefficient of
$e^\l$ equals $1$.
\end{Def}
In fact, the coefficients $a_\mu$ in the expansion
$$
P_\l(q,t)=\sum_{\mu\in\W(\l)} a_\mu E_\mu(q,t)
$$
can be computed explicitly (see, for example, \cite{masterisque}
or \cite[Theorem 3.20, Theorem 4.11]{thesis}). The formulas for
the coefficients $a_\mu$ in the case of a nonreduced root system
are slightly more complicated and we will not list them here. We only
stress that the Proposition \ref{demazure}  is valid for all root
systems.

\begin{Prop}\label{symmcoeff}
If the affine root system $R$ is reduced, then
$$
a_\mu=\prod_{\stackrel {\a\in \RR^+,}
{\scriptscriptstyle{(\a,\mu)>0}}}\frac{t_\a^{-1}-{\bf
q}^{-(\a,\overline\mu)}}{1-{\bf q}^{-(\a,\overline\mu)}}
$$
\end{Prop}

Originally, the definition of symmetric Macdonald polynomials, due
to Macdonald in the reduced setup and to Koornwinder in the
nonreduced setup, preceded that of the nonsymmetric ones. Let us note that when the root system
$R$ is nonreduced the associated polynomials are called in the literature Koornwinder polynomials.

%%%%%%%%%%%%%%%%%%%%%%%%%%%%%%%%%%%%%%%%%%%%%%%%%%%%%%%%%%%%%%%%%%%%%
%%%%%%%%%%%%%%%%%%%%%%%%%%%%%%%%%%%%%%%%%%%%%%%%%%%%%%%%%%%%%%%%%%%%%
%%%%%%%%%%%%%%%%%%%%%%%%%%%%%%%%%%%%%%%%%%%%%%%%%%%%%%%%%%%%%%%%%%%%%
\subsection{The normalization factor}

Given a weight $\l$ let us define the following normalization
factor
$$
e_\l:=\prod_{\stackrel {\a\in R_{nm}^+,}
{\scriptscriptstyle{(\a,\l+\L_0)<0}}} (1-{\bf
q}^{(\a,\overline\l)})
$$
Let us remark that if $\a$ is a positive affine root such that
$(\a,\l+\L_0)<0$ then ${\bf q}^{(\a,\overline\l)}$ is a monomial
in $q^{-1}$, $t_1^{-1}, \cdots, t_{n-1}^{-1}$ and
$(t^*_nt_n)^{-\frac{1}{2}}$. This fact easily follows from the
definition of ${\bf q}^{(\a,\overline\l)}$ and Lemma \ref{lemma5}.
 Therefore, we can state  the following
\begin{Lm}\label{normalizing}
Let $\l$ be a weight. The element $e_\l$ is a polynomial in the
variables $q^{-1}$, $t_1^{-1}, \cdots, t_{n-1}^{-1}$ and
$(t^*_nt_n)^{-\frac{1}{2}}$ with integer coefficients.
\end{Lm}
\begin{Def}
For all weights $\l$ and for all anti--dominant $\mu$, the
polynomials $e_\l E_\l(q,t)$, and respectively $e_\mu P_\mu(q,t)$,
will be called here the normalized  nonsymmetric, respectively
symmetric, Macdonald polynomials.
\end{Def}

These normalized polynomials do not seem to have any particularly interesting properties except for those stated 
in Theorem \ref{polynomiality}.  In type $A$, the normalization factor is much larger 
than the one used to define the integral forms $J_\l(x;q,t)$ of symmetric Macdonald polynomials \cite[VI.$\S 8$]{macsym}, and the same is true for the integral form of the nonsymmetric polynomials \cite[Corrollary 5.2]{kn}.

The following result explains the relationship between the
normalization factors associated to weights in the same
$\W$--orbit.

\begin{Lm}\label{symmnorm}
Let $\l$ be an anti--dominant weight and let $\mu$ be an element
of $\W(\l)$. Then
$$
e_\l=e_\mu \prod_{\stackrel {\a\in \RR_{nm}^+,}
{\scriptscriptstyle{(\a,\mu)>0}}} (1-{\bf q}^{-(\a,\overline\mu)})
$$
\end{Lm}
\begin{proof} The affine root system $R_{nm}$ is  reduced and
hence there will be no loss of generality if we assume $R$ to be
reduced. Let us note first that our hypothesis and Lemma
\ref{lemma4} imply that $w_\l=\w_\mu^{-1}w_\mu$ and
$\ell(w_\l)=\ell(\w_\mu)+\ell(w_\mu)$. Hence, from (\ref{Pi}) we
obtain that
$$
\Pi(w_\l^{-1})=\Pi(\w_\mu)\cup\w_\mu^{-1}(\Pi(w_\mu^{-1}))
$$
Also, the condition $\a\in R^+,~(\a,\l+\L_0)>0$ is equivalent, by
Lemma \ref{lemma2}, to the condition $\a\in \Pi(w_\l^{-1})$.
Therefore,
\begin{eqnarray*}
e_\l&=&\prod_{\a\in \Pi(w_\mu^{-1})} (1-{\bf
q}^{(\w_\mu^{-1}(\a),\overline\l)})\prod_{\a\in \Pi(\w_\mu)}
(1-{\bf q}^{(\a,\overline\l)})\\
&=& \prod_{\a\in \Pi(w_\mu^{-1})} (1-{\bf q}^{(\a,\overline\mu)})
\prod_{\a\in -\w_\mu^{-1}(\Pi(\w_\mu^{-1}))}
(1-{\bf q}^{(\a,\overline\l)})\\
&=& e_\mu \prod_{\a\in \Pi(\w_\mu^{-1})}
(1-{\bf q}^{-(\w_\mu^{-1}(\a),\overline\l)})\\
&=& e_\mu \prod_{\a\in \Pi(\w_\mu^{-1})} (1-{\bf
q}^{-(\a,\overline\mu)})
\end{eqnarray*}
and our statement is proved.
\end{proof}

\begin{Lm}\label{indnorm}
Let $\mu$ be a weight and $s_i$ an affine simple reflection such that $s_i\cdot \mu>\mu$.  Then,
$$
e_{s_i\cdot\mu}=(1-{\bf q}^{-(\a_i^*,\overline\mu)})e_\mu
$$
\end{Lm}
\begin{proof} As above, the affine root system $R_{nm}$ is  reduced and
we can safely assume that $R$ is
reduced. Denote $\l:=s_i\cdot\mu$.
From Lemma \ref{lemma3} we know that $(\mu+\L_0,\a_i)>0$,  $s_iw_\mu=w_\l$ and
$\ell(w_\mu)+1=\ell(w_\l)$. Therefore, if we choose a reduced
decomposition $s_{j_p}\cdots s_{j_1}$ of $w_\mu^{-1}$, then
$s_{j_p}\cdots s_{j_1}s_i$ is a reduced decomposition of
$w_\l^{-1}$ and formula (\ref{Pi}) implies that
$$
\Pi(w_{\l}^{-1})=\{\a_i\}\cup s_i(\Pi(w^{-1}_\mu))
$$
Now, from Lemma \ref{lemma2} we obtain that
$$
\{\a\in R^+~|~ (\a,\l+\L_0)<0\} = \{\a_i\}\cup \{s_i(\a)~|~\a\in
R^+ ~\text{and}~ (\a,\mu+\L_0)<0\}
$$
Therefore,
$$
e_\l=(1- {\bf q}^{(\a_i,\overline\l)})\prod_{\stackrel {\a\in
R^+,} {\scriptscriptstyle{(\a,\mu+\L_0)<0}}} (1-{\bf
q}^{(s_i(\a),\overline\l)})
$$
Note that $(\a_i,\l+\L_0)=-(\a_i,\mu+\L_0)$ and
$(\a_i,\w_\l(\l_j^\vee))=-(\a_i,\w_\mu(\l_j^\vee))$ (we have used
here the first part of Lemma \ref{lemma3}). This implies that
$${\bf q}^{(\a_i,\overline\l)}={\bf
q}^{-(\a_i,\overline\mu)}$$ The same argument shows that
$$\prod_{\stackrel {\a\in R^+,}
{\scriptscriptstyle{(\a,\mu+\L_0)<0}}} (1-{\bf
q}^{(s_i(\a),\overline\l)})=e_\mu$$ which finishes the
proof.
\end{proof}

%%%%%%%%%%%%%%%%%%%%%%%%%%%%%%%%%%%%%%%%%%%%%%%%%%%%%%%%%%%%%%%%%%%%%%%%%%%%%
%%%%%%%%%%%%%%%%%%%%%%%%%%%%%%%%%%%%%%%%%%%%%%%%%%%%%%%%%%%%%%%%%%%%%%%%%%%%%
%%%%%%%%%%%%%%%%%%%%%%%%%%%%%%%%%%%%%%%%%%%%%%%%%%%%%%%%%%%%%%%%%%%%%%%%%%%%%

\subsection{Intertwiners}
The main technical tools used in this paper are the intertwining operators of double affine Hecke algebras defined by Cherednik
\cite{c4}. Our notation and normalization of the intertwiners differ slightly from \cite{c4}, but are consistent with \cite{ion}. The novelty is the
intertwiner denoted below by $\widetilde G_{0,\l}$. 

For any weight $\l$ and any $1\leq i\leq n$ define the operator
$G_{i,\l}= G_{i,\l}(q,t)$ as follows. If  $2\a_i\not \in R$ then
\begin{equation}\label{eq4}
\hspace*{-5.4cm} G_{i,\l}:=(1-{\bf
q}^{-(\a_i,\overline\l)})t_i^{-\frac{1}{2}}T_{i}+ {{\bf
q}^{-(\a_i,\overline\l)}}(1-t_i^{-1})
\end{equation}
If  $2\a_n \in R$ then
\begin{equation}\hspace*{-0.43cm}
G_{n,\l}:= (1-{\bf
q}^{-(\a^*_n,\overline\l)})t_n^{-\frac{1}{2}}T_{n}+ {\bf
q}^{-(\a^*_n,\overline\l)}(1-t_n^{-1})+ {\bf
q}^{-(\a_n,\overline\l)}t_n^{-\frac{1}{2}}(t_{01}^{\frac{1}{2}}-t_{01}^{-\frac{1}{2}})
\end{equation}
The operator $G_{0,\l}$ is defined by the first formula below if
$2\a_0$ is not a root or by the second formula otherwise
\begin{equation}\label{eq5}
\hspace*{-3.1cm} G_{0,\l}:=q^{-(\a_0,\l+\L_0)}\left((1-{\bf
q}^{-(\a_0,\overline\l)})t_{0}^{-\frac{1}{2}}T_{02}+ {{\bf
q}^{-(\a_0,\overline\l)}}(1-t_{0}^{-1})\right)
\end{equation}
\begin{equation}\hspace*{-2cm}
G_{0,\l}:= q^{-(\a_0,\l+\L_0)}\left((1-{\bf
q}^{-(\a^*_0,\overline\l)})t_{01}^{-\frac{1}{2}}T_{02}+ {\bf
q}^{-(\a^*_0,\overline\l)}t_{01}^{-\frac{1}{2}}(t_{02}^{\frac{1}{2}}-t_{02}^{-\frac{1}{2}})\right)
\end{equation}
\begin{equation*}\hspace*{-4cm}
+ q^{-(\a_0,\l+\L_0)}{\bf
q}^{-(\a_0,\overline\l)}t_{01}^{-\frac{1}{2}}(t_{03}^{\frac{1}{2}}-t_{03}^{-\frac{1}{2}})
\end{equation*}
Note that $G_{0,\l}$ differs from the  operator $G_{0,\l}$ in \cite{ion} by a factor of  $q^{-(\a_0,\l+\L_0)}$.
The following result was proved in \cite{ion}.
\begin{Thm}\label{recursion1}
Let $\l$ be a weight and let $\a_i$ be a simple affine root such
that $(\l+\L_0,\a_i)> 0$. Then,
\begin{equation}\label{Grecursion1}
G_{i,\l}\cdot E_\l(q,t)= (1-{\bf
q}^{-(\a_i^*,\overline\l)})E_{s_i\cdot \l}(q,t)
\end{equation}
\end{Thm}

The key role in what follows will be played by the operators
$\widetilde G_{0,\l}$ which are closely related to the operators
$G_{0,\l}$ defined above. The first formula below defines
$\widetilde G_{0,\l}$ in the case of reduced root systems and the
second formula defines it for nonreduced root systems
\begin{equation}\label{eq6}
\hspace*{-3.25cm} \widetilde G_{0,\l}:=({\bf
t}^{(\th,\overline\l)}-{
q}^{-(\a_0,\l+\L_0)})t_{0}^{-\frac{1}{2}}T_{03}+ {{
q}^{-(\a_0,\l+\L_0)}}(1-t_{0}^{-1})
\end{equation}
\begin{equation}\hspace*{0cm}
\widetilde G_{0,\l}:= ({\bf t}^{c_0^{-1}(\th,\overline
\l)}-q^{-(\a_0,\l+\L_0)}{\bf
q}^{-(\a_0,\overline\l)})t_{01}^{-\frac{1}{2}}T_{03}+ {
q}^{-(\a_0,\l+\L_0)}t_{01}^{-\frac{1}{2}}(t_{02}^{\frac{1}{2}}-t_{02}^{-\frac{1}{2}})
\end{equation}
\begin{equation*}\hspace*{-3cm}
+ q^{-(\a_0,\l+\L_0)}{\bf
q}^{-(\a_0,\overline\l)}t_{01}^{-\frac{1}{2}}(t_{03}^{\frac{1}{2}}-t_{03}^{-\frac{1}{2}})
\end{equation*}

\begin{Prop}\label{PropGtildeisG}
Let $\l$ be a weight for which $(\l+\L_0,\a_0)> 0$. Then,
\begin{equation}\label{GtildeisG}
\widetilde G_{0,\l}\cdot E_\l(q,t)=G_{0,\l}\cdot E_\l(q,t)
\end{equation}
\end{Prop}
\begin{proof}
 The formula for the
action of $G_{0,\l}$ involves $T_{02}$, but thanks to Remark
\ref{remarkT02} we can express $T_{02}$ in terms of $T_{03}$ as
follows
$$
T_{02}=T_{02}^{-1}+t_{02}^{\frac{1}{2}}-t_{02}^{-\frac{1}{2}}=T_{03}Y_{\a_0}+
t_{02}^{\frac{1}{2}}-t_{02}^{-\frac{1}{2}}
$$
Note that $T_{03}Y_{\a_0}\cdot E_\l(q,t)={\bf
q}^{(\a_0,\overline\l)}T_{03}\cdot E_{\l}(q,t)$ and therefore
$$
T_{02}\cdot E_\l(q,t)={\bf q}^{(\a_0,\overline\l)}T_{03}\cdot
E_{\l}(q,t)+(t_{02}^{\frac{1}{2}}-t_{02}^{-\frac{1}{2}})E_\l(q,t)
$$
Our claim is an immediate consequence of this formula.
\end{proof}

The normalization of Macdonald polynomials introduced in the previous section is nicely compatible with the action of intertwiners.
\begin{Cor}\label{recursion4}
Let $\l$ be a weight and let $\a_i$ be a simple affine root such
that $s_i\cdot \l>\l$. Then,
\begin{equation}\label{normrecursion}
G_{i,\l}\cdot e_\l E_\l(q,t)= e_{s_i\cdot \l}E_{s_i\cdot \l}(q,t)
\end{equation}
\end{Cor}
\begin{proof}
Straightforward from Theorem \ref{recursion1} and Lemma \ref{indnorm}.
\end{proof}
%%%%%%%%%%%%%%%%%%%%%%%%%%%%%%%%%%%%%%%%%%%%%%%%%%%%%%%%%%%%%%%%%%%%%%%%%%%%%
%%%%%%%%%%%%%%%%%%%%%%%%%%%%%%%%%%%%%%%%%%%%%%%%%%%%%%%%%%%%%%%%%%%%%%%%%%%%%
%%%%%%%%%%%%%%%%%%%%%%%%%%%%%%%%%%%%%%%%%%%%%%%%%%%%%%%%%%%%%%%%%%%%%%%%%%%%%
\subsection{A weak polynomiality result}

Generically, the coefficients of Macdonald polynomials are rational
functions in $q$ and $t$. Therefore, when assigning specific values 
to parameters one has to make sure that the coefficients are well defined for those values. 
One particular instance was considered in \cite[Section 2.3]{ion} 
where it was shown that the limit $E_\l(q,\infty):=\lim_{t\to \infty}E_\l(q,t)$ is well defined
for any weight.  This a fact  is a consequence of the following more precise description of the 
coefficients: they are quotients of polynomials in $q$, $q^{-1}$ and
$t^{-1}$ and the denominators approach 1 when $t\to\infty$.

The technique used to show such a result was introduced  by Knop and Sahi 
\cite{knopsahi} (for Jack polynomials), Knop \cite{kn} (in type $A$) and Cherednik \cite{c4} (general case). The idea is to analyze the action of intertwiners on nonsymmetric polynomials and prove the statement by induction. In the general case, the first result on the  nature of the coefficients of $e_\l E_\l(q,t)$ is due to Cherednik \cite[Corollary 5.3]{c4}: they are Laurent polynomials in $q,t$. A slight improvement (obtain merely by revisiting Cherednik's argument) appeared in 
\cite[Section 2.3]{ion}: the coefficients are polynomials in $q, q^{-1}, t^{-1}$. However, since here we are interested in the limit 
$q\to \infty$ a stronger version of these results is needed.  The general lines of the argument are the same, but the key new ingredient is the operator 
$\widetilde G_{0,\l}$ whose action, unlike that of $G_{0,\l}$, is virtually independent on $q$.

In type $A$, stronger results are know \cite[Corollary 5.2]{kn}, \cite{sahi2}  as the polynomiality of the coefficients is obtained for a normalization of the nonsymmetric polynomials by a significantly smaller factor. The stronger results are obtained along the same lines as below but taking advantage of the additional stability of these polynomials in type $A$. However, for our present purposes the following result will be sufficient.

\begin{Thm}\label{polynomiality}
If the root system $\RR$ is reduced then,

\noi (1) For any weight $\l$ the coefficients of $ e_\l E_\l(q,t)$ are polynomials in  $q^{-1},~t_s^{-1},~t_\ell^{-1}$ with integer coefficients.

\noi (2) For any   anti--dominant weight $\l$, the coefficients of $
e_\l P_\l(q,t)$ are polynomials in $q^{-1},~t_s^{-1},~t_\ell^{-1}$ with integer coefficients.

\noi If the root system $\RR$ is nonreduced, then

\noi (3) For any weight $\l$, the coefficients of $ e_\l E_\l(q,t)$ are polynomials in $~q^{-1}$, $t_{01}^{-1}$, $t^{-1}$, $t_n^{-\frac{1}{2}}t_{01}^{\pm\frac{1}{2}}$,
$t_n^{-\frac{1}{2}}t_{03}^{\pm\frac{1}{2}}$,
$t^{-\frac{1}{2}}_{01}t_{02}^{\pm\frac{1}{2}}$,
$t^{-\frac{1}{2}}_{01}t_{03}^{\pm\frac{1}{2}}$
 with integer
coefficients.
\end{Thm}

\begin{proof}
(1) Let us note first that the proof of part (3) of our statement
follows by precisely the same argument presented below by only
keeping in mind that the elements $T_n$, $T_{01}$, $T_{02}$,
$G_{0,\l}$ and $\widetilde G_{0,\l}$ act in a slightly different
way. In what follows let us assume that $R$ is a reduced root
system.

The statement will be proved by induction on the Bruhat order of
$P$. The minimal elements with respect to the
Bruhat order are the minuscule weights $\l\in \cal O_P$. For such
an element we have that $E_\l(q,t)=e^\l$, by Proposition
\ref{Eminuscule}. Moreover, $\l$ being an element of the affine
fundamental chamber $\cal C$ it satisfies $(\l+\L_0,\a_i)\geq 0$
for all affine simple roots $\a_i$ and in consequence
$(\l+\L_0,\a)\geq 0$ for all affine positive roots. In conclusion,
the normalizing factor $e_\l$ equals 1 and it clear that in this
case $e_\l E_\l(q,t)$ has the predicted properties.

Assume now that $\l$ is an arbitrary non--minuscule weight and
that our statement is true for all weights $\mu<\l$. Since the
weight $\l$ does not belong to the affine fundamental  chamber we
can find an affine simple root $\a_i$ such that
$(\l+\L_0,\a_i)<0$. Let us consider the weight
$$
\mu=s_i\cdot \l
$$
It is clear that $(\mu+\L_0,\a_i)>0$ and therefore Lemma
\ref{lemma3} implies
$$
\l=s_i\cdot\mu >\mu
$$
In particular,  the induction hypothesis applies and we have that
$e_\mu E_\mu(q,t)$ has coefficients which satisfy the conclusion
of the Theorem. Moreover, Corollary \ref{recursion4} implies that
\begin{equation*}
G_{i,\mu}\cdot e_\mu E_\mu(q,t)= e_{\l}E_{\l}(q,t)
\end{equation*}

From the fact that $(1- {\bf q}^{-(\a_i,\overline\mu)})$ appears
as a factor in $e_\l$ and from Lemma \ref{normalizing} we deduce
that ${\bf q}^{-(\a_i,\overline\mu)}$ is a monomial in $q^{-1}$,
$t_s^{-1},~t_\ell^{-1}$. If $i\not =0$, it can be seen directly
from the the formula (\ref{eq4}) that the action of $G_{i,\mu}$
involves only  $t_i^{-1}$ and ${\bf q}^{-(\a_i,\overline\mu)}$. In
conclusion, since the coefficients of $e_\mu E_\mu(q,t)$ are
polynomials in  $q^{-1}$, $t_s^{-1}$, $t_\ell^{-1}$
with integer coefficients, the coefficients of
$e_\l E_\l(q,t)$ have the same property.

If $i=0$, the action of $t_{0}^{-\frac{1}{2}}T_{02}$ involves $q$
in highly nontrivial manner but, nevertheless, it involves only
$t_0^{-1}$, so we can deduce that the coefficients of $e_\l
E_\l(q,t)$ are polynomials in  $q^{\pm 1}$,
$t_s^{-1}$, $t_\ell^{-1}$ with integer
coefficients.

However, from Proposition \ref{GtildeisG}  we also obtain that
$$
\widetilde G_{0,\mu}\cdot e_\mu E_\mu(q,t)=e_\l E_\l(q,t)
$$
From Proposition \ref{HXHY} we deduce that
$T_{03}=X_{\th}T^{-1}_{s_\th}$ and therefore its action involves
the parameters $t_j^{\pm\frac{1}{2}}$ in a complicated way but it
does not involve the parameter $q$ at all. From this fact and from
formula (\ref{eq6}) we obtain that the coefficients of $e_\l
E_\l(q,t)$ are polynomials in  $q^{-1}$,
$t_s^{\pm\frac{1}{2}}$, $t_\ell^{\pm\frac{1}{2}}$ with integer coefficients.

Combining he conclusions of the previous two paragraphs we can conclude that
the coefficients of $e_\l E_\l(q,t)$ must be polynomials in $q^{-1}$, $t_s^{-1}$, $t_\ell^{-1}$
with integer coefficients.

(2) From Proposition \ref{symmcoeff} we deduce that
\begin{equation}\label{symmcoeff2}
e_\l P_\l(q,t)=\sum_{\mu\in\W(\l)}~b_\mu e_\mu E_\mu(q,t)
\end{equation}
with the coefficients
$$b_\mu=a_\mu\frac{e_\l}{e_\mu}=\prod_{\a\in \Pi(\w_\mu^{-1})}
(t_\a^{-1}-{\bf q}^{-(\a,\overline\mu)})$$

Above we have used Lemma \ref{symmnorm} and the definition of
$a_\mu$. We can conclude, using Lemma \ref{lemma2}, that $b_\mu$
is a polynomial in $q^{-1}$, $t_s^{-1}$, $t_\ell^{-1}$ with integer coefficients and therefore by part (1)
we obtain the desired result.
\end{proof}

%%%%%%%%%%%%%%%%%%%%%%%%%%%%%%%%%%%%%%%%%%%%%%%%%%%%%%%%%%%%%%%%%%%%%%%%%%%%%
%%%%%%%%%%%%%%%%%%%%%%%%%%%%%%%%%%%%%%%%%%%%%%%%%%%%%%%%%%%%%%%%%%%%%%%%%%%%%
%%%%%%%%%%%%%%%%%%%%%%%%%%%%%%%%%%%%%%%%%%%%%%%%%%%%%%%%%%%%%%%%%%%%%%%%%%%%%
\section{The $p$--adic degeneration}\label{section3}

\subsection{The limit $q\to\infty$}\label{thelimit}
In recognition of their interpretation within the framework of the representation theory of $p$--adic reductive groups we collect in this section a few results regarding the limit of nonsymmetric Macdonald polynomials as $q\to \infty$.

One important consequence of
Theorem \ref{polynomiality} is that for any weight $\l$ the
coefficients of $E_\l(q,t)$ are rational functions in $q$ and $t$
which can be written as quotients of two polynomials in $q^{-1}$
and $t^{-1}$ (if the root system in question is nonreduced the
statement about the polynomiality in $t^{-1}$ should be altered in
accordance with Theorem \ref{polynomiality}). Moreover, the
denominator $e_\l$ approaches 1 when $q\to \infty$. Therefore, all
coefficients of $E_\l(q,t)$ have finite limits as $q\to \infty$.
Essentially the same argument 
shows that $P_\l(\infty,t)$ is well  defined.

The limit of the nonsymmetric Macdonald polynomials $E_\l(q,t)$ as
the parameter $q$ approaches infinity will be denoted by
$E_\l(\infty,t)$ and will be referred to as  {\sl the $ p$--adic
degeneration} of nonsymmetric Macdonald polynomials. The
terminology is motivated by the fact that for specific values of
the parameters $t_i$ they do have an interpretation as Satake
transforms of some matrix coefficients in unramified principal
series representations of simple $p$--adic groups \cite{ion2}.

The symmetric polynomials $P_\l(\infty,t)$ are in fact already
familiar objects in the representation theory of $p$--adic groups:
up to a scalar factor they are the polynomials that give the
values of zonal spherical functions on a simple  algebraic group
$G$ (defined over a $p$--adic field $k$), relative to a special
maximal compact subgroup $K$, such that the affine root system
associated to $G$ is  the dual affine root system $R^\vee$ or $R$
depending on whether $G$ does or does not split over the
unramified closure of $k$. The parameters $t_i$ represent here
specific integer powers of the cardinality of the residue field of
$k$. If $\RR$ is the root system of type $A_n$ the polynomials
$P_\l(\infty,t)$ are also known as the Hall--Littlewood
polynomials. For reduced root systems, by further specializing the
parameters $t\to \infty$, we obtain that $P_\l(\infty,\infty)$ is
the irreducible Weyl character with lowest weight $\l$ for the
simple complex Lie algebra $\mf g$ with root system $\RR$ (or,
equivalently, for the simple, simply--connected, compact Lie group
with root system $\RR$).

Since we specialized the parameter $q$, from now on we will assume
that the affine Hecke algebra $\H_X$ is defined over $\F_t$ rather
than over $\F$.

\begin{Thm}\label{recursion2}
Let $\l$ be a weight and let $\a_i$ be a simple affine root such
that $(\l+\L_0,\a_i)> 0$. If $i\not = 0$ then
\begin{equation}\label{Trecursion1}
T_{i}\cdot E_\l(\infty,t)= t_i^{\frac{1}{2}}E_{s_i\cdot
\l}(\infty,t)
\end{equation}
If $i=0$ then
\begin{equation}\label{Trecursion2}
T_{03}\cdot E_\l(\infty,t)= t_{01}^{\frac{1}{2}}{\bf
t}^{-c_0^{-1}(\th,\overline\l)}E_{s_0\cdot\l}(\infty,t)
\end{equation}
\end{Thm}
\begin{proof}
The action of the operators $G_{i,\l}$ for $i\not =0$ and
$\widetilde G_{0,\l}$ admit limits as $q\to \infty$. The above
remarks, Theorem \ref{recursion1} and Proposition
\ref{PropGtildeisG} imply the desired result.
\end{proof}

Let $\xi:\H_X \to \F_t$ be the $\Qu$--algebra morphism of
which acts as identity on the parameters $t_i$, $1\leq i\leq n$,
sends $t_{03}$ to $t_{01}$ and
\begin{equation}
\xi(T_i)=t_i^{\frac{1}{2}},\ i\not=0\ \
\text{and}\ \ \xi( T_{03})=t_{01}^{\frac{1}{2}}
\end{equation}
We will abuse notation and write $\xi(w)$ to refer to $\xi(T_w)$
for $w$ in $W$. If we set all the parameters equal $t_1=\cdots=t_n=t_{01}=:t$, then
\begin{equation}\label{eq13}
\xi(w)=t^{{\ell(w)}/{2}}
\end{equation}
Given a weight $\l$ define the following normalization factor
\begin{equation}
f_\l:=\xi(w_\l){\bf t}^{(\l,\overline\l)}
\end{equation}
\begin{Prop}\label{recursion3}
Let $\l$ be a weight and let $\a_i$ be a simple affine root such
that $(\l+\L_0,\a_i)> 0$. Then
$$
T_{i}\cdot f_\l E_\l(\infty,t)= f_{s_i\cdot \l}E_{s_i\cdot
\l}(\infty,t)
$$

\end{Prop}
\begin{proof}
From Theorem \ref{recursion2} it is clear that we only have to
check that under our hypothesis
$f_{s_i\cdot\l}=t_i^{\frac{1}{2}}f_\l$ if $i\not = 0$ and
$f_{s_0\cdot\l}=t_{01}^{\frac{1}{2}}{\bf
t}^{-c_0^{-1}(\th,\overline\l)}f_\l$ if $i=0$.

Let us assume first that $i\not =0$. Then, Lemma \ref{lemma3}
implies that $\xi(w_{s_i\cdot\l})=t_i^{\frac{1}{2}}\xi(w_\l)$ and
$\w_{s_i\cdot\l}=s_i\w_\l$ which prove the above claim. If $i=0$
then by the same result
$\xi(w_{s_0\cdot\l})=t_{01}^{\frac{1}{2}}\xi(w_\l)$ and
$\w_{s_0\cdot\l}=s_\th\w_\l$. In particular, for any $x\in\Gc^*$
\begin{eqnarray*}
(s_0\cdot\l, \w_{s_0\cdot\l}(x))&=&(s_\th(\l)+c_0^{-1}\th,
s_\th\w_{\l}(x))\\
&=&(\l,\w_{\l}(x))-(c_0^{-1}\th, \w_{\l}(x))
\end{eqnarray*}
Therefore, ${\bf t}^{(s_0\cdot\l,\overline{s_0\cdot\l})}={\bf
t}^{-c_0^{-1}(\th,\overline\l)}{\bf t}^{(\l,\overline\l)}$ and the
proof is completed.
\end{proof}
\begin{Cor}\label{cor4}
Let $\l$ be a weight. Then
$$
T_{w_\l}\cdot f_{\tilde \l}e^{\tilde \l}=f_\l E_\l(\infty,t)
$$
\end{Cor}
\begin{Cor}\label{cor2}
For any weight $\l$, the parameters $t_i$, $t_{01}$
appear in $\xi({w_\l})f_\l/f_{\tilde \l}$ with integer exponents. If all the
parameters are equal (and denoted by $t$) then
$$
f_\l=t^{-\ell(\w_\l)/2}
$$
\end{Cor}
\begin{proof}
The first claim follows immediately from the fact that the
monomials $f_\mu$ are obtained inductively as in the proof of the
above Proposition and the fact that
$c_0^{-1}(\th,\w_\mu(\l_i^\vee))$ are integers.

For the second claim, note that by using the first part of Lemma
\ref{lemma3} and formula (\ref{lusztiglength}) we obtain that
$${\bf t}^{(\l,\overline\l)}={\bf t}^{(\l_-,\overline{\l_-})}=t^{-\ell(\tau_{\l_-})/2}$$
Now, using the third part of Lemma \ref{length} and the third part
of Lemma \ref{lemma4} we deduce that
$$
{\bf
t}^{(\l,\overline\l)}=t^{-\frac{1}{2}(\ell({w_\l})+\ell(\w_\l))}
$$
and our statement immediately follows.
\end{proof}

 The next result expresses the
relationship between the $p$--adic degeneration of the symmetric
and nonsymmetric Macdonald polynomials.
\begin{Prop}\label{demazure}
Let $\l$ be an anti--dominant weight. Then,
$$
P_\l(\infty,t)=\sum_{\mu\in\W(\l)}\xi(\w_\mu)^{-2}
E_\mu(\infty,t)
$$
\end{Prop}
\begin{proof} Straightforward from Proposition \ref{symmcoeff} and
the definition of $\xi(\w_\mu)$. We only note that although we
only stated Proposition \ref{symmcoeff} for reduced root systems,
a similar fact holds for nonreduced root systems \cite[Theorem
4.11]{thesis} and the limit of the corresponding coefficients
$a_\mu$ in the limit $q\to\infty$  equals also
$\xi(\w_\mu)^{-2}$.
\end{proof}

\begin{Prop}\label{thm3}
Assume $\l$ is a dominant weight. Then, $E_\l(\infty,t)=e^\l$.
\end{Prop}
\begin{proof}  First, note that Corollary \ref{cor4} allows
us to write
\begin{equation}
T_{w_\l}\cdot f_{\tilde \l} e^{\tilde \l}=f_\l E_\l(\infty,t)
\end{equation}
Second, from Lemma \ref{lemma4} we know that
$\ell(\tau_\l)=\ell(w_\l)+\ell(\w_\l^{-1})$. Hence,
$T_{\tau_\l}=T_{w_\l}\omega_{\tilde \l}T_{\w_\l^{-1}}$. The weight
$\l$ being dominant we have that $X_\l=T_{\tau_\l}$ and
consequently,
$$
T_{w_\l}X_{\tilde \l}=X_\l
T_{\w_\l^{-1}}^{-1}T_{\w_{\tilde\l}^{-1}}
$$
Therefore, $$T_{w_\l}\cdot f_{\tilde \l} e^{\tilde \l}= f_{\tilde \l} X_\l
T_{\w_\l^{-1}}^{-1}T_{\w_{\tilde\l}^{-1}} \cdot 1 =
f_{\tilde \l}\xi(\w_\l)^{-1}\xi(\w_{\tilde\l})e^\l$$ Also, the
coefficient of $e^\l$ in $E_\l(\infty,t)$ is 1 and the
claim follows.
\end{proof}

One immediate consequence of the above computation is that  
\begin{equation}\label{eq17}
f_\l/f_{\tilde \l}=\xi(\w_\l)^{-1}/\xi(\w_{\tilde\l})^{-1}
\end{equation}
for $\l$ anti-dominant, but keeping in mind (see the proof of Corollary \ref{cor2}) that ${\bf t}^{(\mu,\overline\mu)}$ is constant for $\mu\in \W(\l)$
and the Lemma \ref{lemma4} (1) we deduce that (\ref{eq17}) is true for any weight $\l$.

\begin{Cor}\label{cor1}
Let $\l$ be a weight. Then,
$$
E_\l(\infty,t)=\xi(\w_\l w_\circ )^{-1}T_{\w_\l w_\circ}\cdot e^{\l_+}
$$
In particular, if $\RR$ nonreduced the polynomials $E_\l(\infty, t)$ are free of the variables $t_{01}, ~t_{02}$.
\end{Cor} 

As anticipated in \cite{ion} we have the following
\begin{Cor}\label{cor7}
Assume $\RR$ is reduced. For any  weight $\l$
$$E_\l(\infty,\infty):=\lim_{t\to\infty}E_\l(\infty,t)$$ is the
Demazure character associated to the irreducible representation of
$\mf g$ with highest weight $\l_+$ and extremal weight $\l$. In
particular, if $\l$ is anti--dominant, then $E_\l(\infty,\infty)$
is the Weyl character of the irreducible representation of $\mf g$
with lowest weight $\l$. Moreover,
\begin{equation}\label{eq7}
E_\l(\infty,\infty)=\lim_{t\to\infty}\lim_{q\to\infty}E_\l(q,t)=\lim_{q\to\infty}\lim_{t\to\infty}E_\l(q,t)
\end{equation}
\end{Cor}
\begin{proof}
The fact that $E_\l(\infty,\infty)$ is the Demazure character
associated to the irreducible representation of $\mf g$ with
highest weight $\l_+$ and extremal weight $\l$ is an immediate
consequence of Corollary \ref{cor1} and of the Demazure character
formula \cite{demazure}. The equation (\ref{eq7}) follows by combining this result
with Theorem 3 in \cite{ion}.
\end{proof}
It is also clear that $E_\l(\infty,1)=e^\l$ for all weights $\l$
and therefore the polynomials $E_\l(\infty,t)$ 
interpolate between monomials $e^\l$ and Demazure characters
associated to the irreducible representation of $\mf g$ with
highest weight $\l_+$ and extremal weight $\l$. This property is
the nonsymmetric analogue of the corresponding fact regarding the
symmetric polynomial $P_\l(\infty,t)$ which is know to
interpolate between the symmetrized monomial
$\sum_{\mu\in\W(\l)}e^\mu$ and  the Weyl character of the
irreducible representation of $\mf g$ with lowest weight $\l$
(note that here $\l$ is an anti--dominant weight).
\begin{Cor}\label{cor8}
Let $\l$ be an anti--dominant weight. Then,
\begin{equation}\label{eq8}
P_\l(\infty,t)=\xi(w_\circ)^{-1}\sum_{\mu\in\W(\l)}\xi(\w_\mu)^{-1}T_{\w_\mu w_\circ}\cdot e^{\l_+}
\end{equation}
\end{Cor}
\begin{proof}
Straightforward from Proposition \ref{demazure} and Corollary
\ref{cor1}.
\end{proof}
Remark that if $\RR$ is a reduced root system and $t\to \infty$
the equation (\ref{eq8}) becomes precisely the Demazure character
formula the irreducible representation of $\mf g$ with highest
weight $\l_+$. In the light of the connection between
$P_\l(\infty,t)$ and spherical functions on simple groups over
$p$--adic fields, the equation (\ref{eq8}) could be seen as a
counterpart of Demazure's formula for this type of spherical
functions. The above result also follows from equation (5.4) and
Lemma 4.2 in \cite{knop2} together with Macdonald's formula for
the Satake transforms of the elements $N_\l$ in \cite{knop2}.

It is natural to introduce the following normalization of the nonsymmetric Macdonald polynomials. 
\begin{Def} For any weight $\l$ define $$
{\widetilde E_\l(q,t)}:=\xi(\w_\l)^{-1}E_\l(q,t)
$$
\end{Def}
We close this section with a reformulation of Corollary \ref{cor4} in terms of the above normalization. For roots systems of type $A$ the result was proved by Knop \cite[Corollary 5.3]{kn}.
\begin{Cor}\label{cor9}
Let $\l$ be a weight. Then
$$
T_{w_\l}\cdot {\widetilde E}_{\tilde \l}(\infty,t)={\widetilde E}_\l(\infty,t)
$$
\end{Cor}
\begin{proof}
Clearly, 
\begin{eqnarray*}
T_{w_\l}\cdot {\widetilde E}_{\tilde \l}(\infty,t)&=&({\xi(\w_{\tilde \l})^{-1}}/{f_{\tilde \l}})T_{w_\l}\cdot f_{\tilde \l} {e}^{\tilde \l}\\
&=&({\xi(\w_{\tilde \l})^{-1}}/{f_{\tilde \l}}) f_{\l} E_\l(\infty,t)\\
&=&{\widetilde E}_\l(\infty,t)
\end{eqnarray*}
In  the last step (\ref{eq17}) was used.
\end{proof}

%%%%%%%%%%%%%%%%%%%%%%%%%%%%%%%%%%%%%%%%%%%%%%%%%%%%%%%%%%%%%%%%%%%%%%%%%%%%%
%%%%%%%%%%%%%%%%%%%%%%%%%%%%%%%%%%%%%%%%%%%%%%%%%%%%%%%%%%%%%%%%%%%%%%%%%%%%%
%%%%%%%%%%%%%%%%%%%%%%%%%%%%%%%%%%%%%%%%%%%%%%%%%%%%%%%%%%%%%%%%%%%%%%%%%%%%%
\subsection{Normalized intertwiners}\label{normalized}

For any weight $\l$ and any $0\leq i\leq n$ define the following normalized versions of the intertwiners.
The second formula below defines
$I_{0,\l}$ for reduced root systems and the
third formula defines it for nonreduced root systems
\begin{eqnarray*}
I_{i,\l}&:=&t_i^{\frac{1}{2}}G_{i,\l}/(1-{\bf
q}^{-(\a^*_i,\overline\l)}), \quad i\neq 0\\
 I_{0,\l}&:=&t_{0}^{\frac{1}{2}}\widetilde G_{0,\l}/({\bf
t}^{(\th,\overline\l)}-{
q}^{-(\a_0,\l+\L_0)})  \\
 I_{0,\l}&:=& t_{01}^{\frac{1}{2}}\widetilde G_{0,\l}/({\bf t}^{c_0^{-1}(\th,\overline
\l)}-q^{-(\a_0,\l+\L_0)}{\bf
q}^{-(\a_0,\overline\l)})  
\end{eqnarray*}

Fix  a weight $\l$ and a reduced decomposition $s_{j_\ell}\cdots s_{j_1}$ of $w_\l$. Denote $\l_{(1)}=\tilde \l$ and 
 $\l_{(i)}=s_{j_{i-1}}\cdots s_{j_1}\cdot \tilde \l$ for $2\leq i\leq \ell$. With this notation define
 $$
 I_{w_\l}:=I_{j_\ell,\l_{(\ell)}}\cdots I_{j_1,\l_{(1)}}
 $$
 \begin{Thm}\label{normalized_intertwiners}
Let $\l$ be a weight. Then
$$
I_{w_\l}\cdot {\widetilde E}_{\tilde \l}(q,t)={\widetilde E}_\l(q,t)
$$
\end{Thm}
\begin{proof}   First,  remark that Theorem \ref{recursion2} holds for $I_i$ replacing $T_i$ and $E_\l(q,t)$ replacing $E_\l(\infty,t)$ which then implies the conclusions of Proposition \ref{recursion3} and Corollary \ref{cor4} (under the same substitutions).
The proof can be concluded following exactly the same line as the proof of Corollary \ref{cor9}.
\end{proof}

%%%%%%%%%%%%%%%%%%%%%%%%%%%%%%%%%%%%%%%%%%%%%%%%%%%%%%%%%%%%%%%%%%%%%%%%%%%%%
%%%%%%%%%%%%%%%%%%%%%%%%%%%%%%%%%%%%%%%%%%%%%%%%%%%%%%%%%%%%%%%%%%%%%%%%%%%%%
%%%%%%%%%%%%%%%%%%%%%%%%%%%%%%%%%%%%%%%%%%%%%%%%%%%%%%%%%%%%%%%%%%%%%%%%%%%%%

\section{Bases for maximal parabolic modules}\label{section4}

\subsection{The Kazhdan-Lusztig involution}  In this section we begin to explore the connection between the nonsymmetric Macdonald polynomials and the Kazhdan--Lusztig theory (in its parabolic version \cite{deodhar}). 
We start by recalling the construction of the Kazhdan--Lusztig basis \cite{kl1} in its multi-parameter version \cite{lusztig3} and some other basic facts.

From now on we will work over the field obtained from $\F$ by specializing the parameters $t_{01}, ~t_{02}$ (if present) to 1. As stated in the Corollary \ref{cor1}
the $p$--adic limit is independent of these variables hence unaffected by the specialization. To avoid introduction new notation we will use the old notation for the fields and polynomials under consideration. From the point of view of the Kazhdan--Lusztig theory a new feature is the introduction of $q$ as a parameter. As we will see below, this is completely harmless in regard to the general theory, but it will allow us to draw some conclusions  regarding the interpretation of the nonsymmetric polynomials within this framework.

Let  $
\chi:\H_X\to \F
$
be the $\F$--algebra map  which sends each of the generators $T_i$, $T_{03}$ to the square root of the corresponding parameter. 

 The Kazhdan--Lusztig involution $\kappa$ is the involution of the  algebra
 $\H_X$ which inverts the parameters $q$, $\{t_\a\}_{\a\in R}$ and the generators
 $T_i$, $T_{03}$. On a standard basis element it acts as follows
 \begin{equation}\label{eq18}
 \kappa(T_w)=T_{w^{-1}}^{-1}
 \end{equation}
In fact, we can extend $\kappa$ to $\H_X^e$ (as an algebra map) by letting it act as identity on $\Omega$. The formula (\ref{eq18}) 
is then valid for any $w\in W^e$.

Recall from  \cite[Proposition 2]{lusztig3},  \cite[Theorem 5.2]{lusztig4} the following result.
\begin{Thm}\label{existence}
For any element $w$ of the affine Weyl group $W$ there is a unique
element $C'_w$ of $\H_X$ which satisfies the properties
\begin{enumerate}
\item[(a)] $\kappa(C_w')=C_w'$

\item[(b)] $C_w'=\sum_{y\leq w}
P^*_{y,w}(t)T_y$, where $P^*_{w,w}(t)=1$ and, if $y<w$,
 $P^*_{y,w}(t)$ are polynomials in $\{t^{-\frac{1}{2}}_\a\}_{\a\in R}$ with integer
coefficients and no constant term.
\end{enumerate}
Moreover, $P_{y,w}(t):=\chi(y)^{-1}\chi(w)P^*_{y,w}(t)\in \Z[ t_\a|\a\in R ]$.
\end{Thm}
Lusztig's result is in fact valid for any Coxeter group. 
The polynomials $P_{y,w}(t)$ are  {\sl
Kazhdan--Lusztig polynomials} (for the affine Weyl group $W$). 
For equal parameters the polynomials $P_{y,w}(t)$ have non--negative
coefficients. This fact  follows from a beautiful
cohomological interpretation \cite[Theorem 5.5]{kl2} in terms of
the Deligne--Goresky--MacPherson middle intersection cohomology.

We also need the following basic facts \cite[(4.2), (4.3)]{lusztig3}. 
\begin{Prop}\label{symmetry} 
Let $s_i$ be a simple reflection and $w$ be an element of $W$ such that $s_iw<w$. Then,
$$
(T_i-t_i^{\frac{1}{2}})C'_w=0
$$
\end{Prop}
\begin{Lm}
Let $s_i$ be an affine simple reflection and let $x$, $y$ let
elements of $W$ such that $x<y$, $x<xs_i$ and $ys_i<y$. Then,
$$
P^*_{x,y}(t)=t_i^{-\frac{1}{2}}P^*_{xs_i,y}(t)
$$
\end{Lm}
\begin{Cor}\label{cor3}
Let $\mu$ and $\l$ be two weights such that $\mu\leq\l$. Then,
$$
P^*_{v_\mu y,v_\l}(t)=\chi(y)^{-1}P^*_{v_\mu,v_\l}(t)
$$
for all $y$ in $W_{\tilde \l}$.
\end{Cor}
The elements $C_{v_\l}'$ can be factorized as follows.

\begin{Lm}\label{lemma7}
Let $\l$ be a weight. Then
$$
C_{v_\l}'=\left(\sum_{\mu\leq \l}
P^*_{v_\mu,v_\l}(t)T_{w_\mu}\omega_{\tilde\l}\right)\left(\chi(w_{\circ})^{-1}
\sum_{x\in \W}{\chi(x)}T_x\right)\omega_{\tilde \l}^{-1}
$$
\end{Lm}
\begin{proof}
Let us  remark first that from Lemma \ref{lemma6}
$$
\{y\in W~|~y\leq v_\l\}=\bigcup_{\mu\leq\l} v_\mu W_{\tilde\l}
$$
and from the above Corollary
$$ P^*_{v_\mu y,v_\l}(t)=\chi(y)^{-1}P^*_{v_\mu,v_\l}(t)
$$
for all $y$ in $W_{\tilde \l}$ and $\mu\leq \l$. Hence, the element $C_{v_\l}'$ takes the form
$$
C_{v_\l}'=\sum_{\mu\leq \l}
P^*_{v_\mu,v_\l}(t)\left( \sum_{y\in 
W_{\tilde\l}}{\chi(y)^{-1}}T_{v_\mu y}\right)
$$
Since $v_\mu W_{\tilde\mu}=w_\mu W_{\tilde\mu}$ and $\ell(v_\mu
y)=\ell(v_\mu)-\ell(y)=\ell(w_\mu)+\ell(w_{\circ,\tilde \l})-\ell(y)$ for any $y$ in $ W_{\tilde\mu}$ we get that
$$
\sum_{y\in 
W_{\tilde\l}}{\chi(y)^{-1}}T_{v_\mu y}={\chi(w_{\circ,\tilde \l})^{-1}}T_{w_\mu}
 \sum_{x\in W_{\tilde\l}}{\chi(x)}T_x
$$
Now, $W_{\tilde \l}=\omega_{\tilde \l}\W\omega_{\tilde \l}^{-1}$ and $\chi$ is invariant under the conjugation action of $\Omega$. Our claim now immediately
follows.
\end{proof}

 For the following result recall the notation in Section \ref{normalized} and that we work under the assumption $t_{01}=t_{02}=1$.
 \begin{Lm}\label{lemma8} 
 Let $\l$ be a weight. Then $I_{w_\l}\omega_{\tilde \l}$ is fixed by $\kappa$.
 \end{Lm}
\begin{proof}
It is a straightforward check that each factor of $I_{w_\l}$ is  fixed by $\kappa$.
\end{proof}
%%%%%%%%%%%%%%%%%%%%%%%%%%%%%%%%%%%%%%%%%%%%%%%%%%%%%%%%%%%%%%%%%%%%%%%%%%%%%
%%%%%%%%%%%%%%%%%%%%%%%%%%%%%%%%%%%%%%%%%%%%%%%%%%%%%%%%%%%%%%%%%%%%%%%%%%%%%
%%%%%%%%%%%%%%%%%%%%%%%%%%%%%%%%%%%%%%%%%%%%%%%%%%%%%%%%%%%%%%%%%%%%%%%%%%%%%

\subsection{The parabolic module}  
Restricting $\chi$ to $\h$ we obtain $(\chi_{|\h},\F_t)$ a one dimensional representation of $\h$. The induced representation
$$
\text{ind}_{\h}^{\H_X}(\chi):=\H_X\otimes_{\h}\F_t 
$$
is a left module for $\H_X$. In general, there are several standard maximal parabolic subgroups of $W$ isomorphic to $\W$ (as many as the order of $\Omega$) and one can construct in the same manner the corresponding induced representation of the affine Hecke algebra $\H_X$. These, however, are all isomorphic to the one defined above. One can consider all of them together by constructing the $\H_X^e$--module 
$$
\text{ind}_{\h}^{\H_X^e}(\chi):=\H_X\otimes_{\h}\F_t 
$$
Following Knop \cite{knop} we call the above module the (maximal) parabolic module of $\H_X^e$. By Proposition \ref{HXHY} the parabolic module has a basis given by $\{X_\l\otimes 1\}_{\l\in P}$ and it is isomorphic as a $\H_X^e$--module to $\cal R_t$ ($X_\l\otimes 1$ and $e^\l$ correspond under the isomorphism).

It is a standard fact (see \cite{deodhar}) that there exists an involution (still called the Kazhdan--Lusztig involution and denoted by $\kappa$)
$$
\kappa:\cal R \to \cal R, \quad f\mapsto f^\kappa
$$
compatible with the one on $\H_X^e$ in the following sense
\begin{equation}\label{eq19}
(H\cdot f)^\kappa=\kappa(H)\cdot f^\kappa
\end{equation}
for any $H\in \H_X^e$ and $f\in \cal R$. In our case, however, everything can be made quite explicit. Note that for $\l$ dominant $X_\l=T_{\tau_\l}$ and therefore
\begin{eqnarray*}
\kappa(X_\l)&=&T_{\tau_{-\l}}^{-1}\\
&=&T_{w_\circ}T_{\tau_{-w_\circ(\l)}}^{-1}T_{w_\circ}^{-1}\\
&=&T_{w_\circ}X_{{w_\circ(\l)}}T_{w_\circ}^{-1}
\end{eqnarray*}
The map $\kappa$ being an algebra morphism we obtain that $$\kappa(X_\l)=T_{w_\circ}X_{{w_\circ(\l)}}T_{w_\circ}^{-1}$$ for any weight $\l$. Applying now
(\ref{eq19}) for $H=X_\l$ and $f=1$ we obtain that 
$$
\kappa(e^\l)=\chi(w_\circ)^{-1}T_{w_\circ}\cdot e^{w_\circ(\l)}
$$
Of course, $\kappa$ acts on $\F$ by inverting the parameters.

Keeping in mind that $\{w_\l\omega_{\tilde \l}\}_{\l\in P}$ is the set of minimal coset representatives for $W^e/\W$ we obtain the following bases of $\cal R$ which are induced from elements of $\H^e_X$.
\begin{Def} Define the following are bases for $\cal R$:
\begin{enumerate}
\item[(a)] the standard basis: $\{T_{w_\l}{\omega_{\tilde \l}}\cdot 1\}_{\l\in P}$;
\item[(b)] the dual standard basis: $\{T^{-1}_{w^{-1}_\l}{\omega_{\tilde \l}}\cdot 1\}_{\l\in P}$;
\item[(c)] the canonical basis $\{C'_\l:=t^{\ell(w_\circ)/2}\W(t)^{-1} C'_{v_\l}\omega_{\tilde \l}\cdot 1\}_{\l\in P}$.
\end{enumerate}
Above we denoted by $\W(t):=\sum_{\w\in\W}t^{\ell(\w)}$ the
Poincar\' e polynomial of $\W$.
\end{Def}

The coefficients of the expansion of the canonical basis in the standard basis are the parabolic Kazhdan--Lusztig polynomials (for the maximal parabolics $W_{\tilde \l}\subset W$) of Deodhar \cite{deodhar}.
We are now ready to establish one connection between the nonsymmetric Macdonald polynomials and the (parabolic) Kazhdan--Lusztig theory.

\begin{Thm}\label{state2}  The basis  $\{ {\widetilde E}_\l(q,t)\}_{\l\in P}$ of the parabolic module of the affine Hecke algebra $\H_X^e$ is 
 invariant under the Kazhdan--Lusztig involution. Moreover, 
\begin{enumerate}
\item $\{ {\widetilde E}_\l(\infty,t)\}_{\l\in P}$ is
the standard basis;

\item  $\{ {\widetilde E}_\l(0,t)\}_{\l\in P}$ is
the dual standard basis.
\end{enumerate}
\end{Thm}
\begin{proof} Let $\l$ be a weight and fix a reduced decomposition of $w_\l$.
By Lemma \ref{lemma8} the elements $I_{w_\l}\omega_{\tilde \l}$ are fixed by $\kappa$ and therefore, using (\ref{eq19}) for $H=I_{w_\l}\omega_{\tilde \l}$ and $f=1$, we obtain that 
\begin{eqnarray*}
I_{w_\l}\omega_{\tilde \l}\cdot 1 &=& I_{w_\l}\cdot \widetilde E_{\tilde \l}(q,t)\\
&=& \widetilde E_{\l}(q,t)
\end{eqnarray*}
is fixed by $\kappa$. Also, (1) is exactly Corollary \ref{cor9}. To explain (2) note that
$$
 \widetilde E_{ \l}(q,t)=\chi(w_\circ)^{-1}T_{w_\circ}\cdot w_\circ(\widetilde E_\l(q^{-1},t^{-1}))\\
$$
As the limit as $q\to 0$ of the right hand side exists (see the discussion at the beginning of Section \ref{thelimit}) the limit of the left hand side also 
exists and $$ \widetilde E_{ \l}(0,t)= \widetilde E^\kappa_{ \l}(\infty,t) $$ In conclusion, $\widetilde E_{ \l}(0,t)$ is an element of the dual standard basis.
\end{proof}

The expansion of the canonical basis in terms of the standard basis takes the following form.
\begin{Prop}\label{prop3} Let $\l$ be a weight. Then
\begin{equation}\label{eq20}
C_{\l}'=\sum_{\mu\leq \l}
P^*_{v_\mu,v_\l}(t)\widetilde E_\mu(\infty,t)
\end{equation}
\end{Prop}
\begin{proof}
Staightforward from Lemma \ref{lemma7}.
\end{proof}

It is clear from (\ref{eq20}) that the coefficients of $C'_\l$ are polynomials in $\{t^{-\frac{1}{2}}_\a\}_{\a\in R}$ with integer
coefficients. For reduced root systems and $\l$ anti-dominant $C'_\l$ were shown (originally in  \cite{lusztig2}, later reproved by several authors ) to be Weyl characters of the irreducible representation of $\mf g$ with lowest weight
$\l$.  For completeness, we also give a proof here. The argument follows the idea used in \cite{kato}.
\begin{Thm}\label{KLantidominant}

 Assume $\RR$ to be a reduced root system and let $\l$ be an anti-dominant weight. Then, $C'_\l$ is the Weyl character
of the irreducible representation of $\mf g$ with lowest weight
$\l$.
\end{Thm}
\begin{proof} Let $s_i$ be a simple reflection.
The condition $v_\l>s_iv_\l$ in Proposition \ref{symmetry} translates into $\l\ge s_i\cdot \l$ and this certainly holds for any $1\le i\le n$ (an anti-dominant weight is the highest in its $\W$ orbit). Therefore, 
$$
(T_i-t_i^{\frac{1}{2}})C'_\l=0, \quad \text{for all~~} 1\le i\le n
$$
which is equivalent to $C'_\l$ being $\W$--invariant. On $\W$--invariant elements, $\kappa$ has only the effect of inverting the parameters. Thus, $C'_\l$ being fixed by $\kappa$ and with coefficients  polynomials in $\{t^{-\frac{1}{2}}_\a\}_{\a\in R}$ forces it to be free of parameters. Therefore, we may safely take all the parameters to infinity in (\ref{eq20}) without altering $C'_\l$. The only term on the right hand side which survives this process in $\widetilde E_\l(\infty,\infty)$ which by Corollary \ref{cor7} is the specified Weyl character.
\end{proof}

Explicit formulas or representation--theoretical interpretations of the elements $C'_\l$ beyond the case described above seem to be unknown.  However, in the equal parameter case a few special properties are expected.

For the following remarks assume that  $\RR$ is a reduced root system and the parameters are equal $t_s=t_\ell=:t$. Computational evidence suggests the following
\begin{conjecture}\label{conjecture}
 For any weight $\l$ the polynomial $C'_\l$ is the $T$--character of a graded $B$--module. In particular, the positive integers $P_{v_\mu,v_\l}(1)$ represent weight multiplicities in $B$--modules.
\end{conjecture}

Given that the  the polynomials $\widetilde E_\l(q,t)$ interpolate between the standard and the dual standard basis the expansion of the canonical basis in 
terms of them is especially intriguing. One case in particular draws attention: the expansion of the canonical basis (for $\RR$ reduced, equal parameters)
in the polynomials $\widetilde E_\l(q,t)$  (for $q=t, ~t_s=t, ~t_\ell=t^r$) seems to characterized by a support condition which in turn reduces the problem of computing all the elements of the canonical basis (infinitely many) for a fixed root system to a finite computation. We will report on these investigations elsewhere. 

%%%%%%%%%%%%%%%%%%%%%%%%%%%%%%%%%%%%%%%%%%%%%%%%%%%%%%%%%%%%%%%%%%%%%%%%%%%%%
%%%%%%%%%%%%%%%%%%%%%%%%%%%%%%%%%%%%%%%%%%%%%%%%%%%%%%%%%%%%%%%%%%%%%%%%%%%%%
%%%%%%%%%%%%%%%%%%%%%%%%%%%%%%%%%%%%%%%%%%%%%%%%%%%%%%%%%%%%%%%%%%%%%%%%%%%%%

\subsection{Orthogonality} We now revert back to the multi-parameter situation. Since the polynomials $E_\l(q,t)$
form a basis of $\cal R$ orthogonal with respect to the scalar
product $\<\cdot,\cdot\>_{q,t}$. It is natural to ask if such a
property holds for the polynomials $E_\l(\infty,t)$ with respect
to the space $\cal R_t$ and a suitable degeneration of the scalar
product $\<\cdot,\cdot\>_{q,t}$ as $q\to\infty$.  
Unfortunately,
the definition of the Cherednik scalar product involves the
involution $\overline{~\cdot~}$ on $\cal R$, which inverts the
parameter $q$ and it is therefore inconsistent with the process of
taking the limit $q\to \infty$. However, we can try to examine the
limit as $q\to\infty$ of
\begin{equation}\label{eq9}
\<E_\l(q,t),E_\mu(q,t)\>_{q,t}=CT\left(E_\l(q,t)\overline{E_\mu(q,t)}C(q,t)\right)
\end{equation}
Although it is clear that the limit as $q$ approaches infinity of
$E_\l(q,t)$ and $C(q,t)$ exists (and equals $E_\l(\infty,t)$ and,
respectively, $C(\infty,t)$) it is not clear what happens to $\overline{E_\mu(q,t)}$ in the limit. Before stating a result of
Cherednik which will allow us to perform such a computation we
need to introduce some notation.

Let $\varsigma$ be the involution of $\cal R$ which fixes the
parameters $q$ and $t$ and, for any weight $\l$, sends $e^\l$ to
$e^{-w_\circ(\l)}$. Also, let
$$\iota:\cal R\to\cal R, \quad \iota=\chi(w_\circ)T_{w_\circ}^{-1}
\varsigma$$ We will also use the notation $f^\iota:=\iota(f)$ for
any element $f$ of $\cal R$.
\begin{Prop}
Let $\l$ be a weight. Then,
\begin{enumerate}
\item
$\overline{E_\l(q,t)}=\chi(\w_\l)^{-2}\chi(w_\circ)T_{w_\circ}^{-1}\cdot
E_{-w_\circ(\l)}(q,t)$

\item $\varsigma(E_\l(q,t))=E_{-w_\circ(\l)}(q,t)$

\item $\overline{\widetilde E_\l(q,t)}={\widetilde E}_\l^\iota(q,t)$
\end{enumerate}
\end{Prop}

The first claim was proved by Cherednik \cite[Proposition 3.3]{c2}
and the second claim follows along exactly the same lines.
Although it is not explicitly stated in \cite{c2} it is implicitly used at
several places. The third claim is simply a combination of the
previous two. As an immediate consequence we have the following
\begin{Cor}
For any weight $\l$, the limit of $\overline{E_\l(q,t)}$ as $q$
approaches infinity exists and equals
$\chi(\w_\l)^{-2}E_\l^\iota(\infty,t)$.
\end{Cor}
Using this result, it is clear that we can take the limit
$q\to\infty$ directly on the right hand side of equation
(\ref{eq9}) and obtain that
\begin{equation}\label{eq10}
\lim_{q\to\infty}\<E_\l(q,t),E_\l(q,t)\>_{q,t}=\chi(\w_\l)^{-2}CT\left(
E_\l(\infty,t)E_\l^\iota(\infty,t)C(\infty,t)\right)
\end{equation}
On the other hand, $E_\l(q,t)$ are orthogonal with respect to the
scalar product $\<\cdot,\cdot\>_{q,t}$ and their norms are
explicitly known (see, for example, \cite{c3} for formulas for
reduced root systems or \cite{macbook} for completely general
results). By examining these norms it is easy to see that
$$
\lim_{q\to\infty}\<E_\l(q,t),E_\l(q,t)\>_{q,t}=1
$$
It is therefore natural to define, following \cite[Corollary
4.3]{c2}, the following symmetric {scalar product} on $\cal R_t$. For $f$
and $g$ in $\cal R_t$ let
\begin{equation}\label{eq11}
\<f,g\>_t:=CT\left( fg^\iota C(\infty,t)\right)
\end{equation}

We have proved the following:
\begin{Prop}\label{prop2}
The polynomials ${\widetilde E_\l(\infty,t)}$ form a basis of
$\cal R_t$ which is orthonormal with respect to the scalar product
$\<\cdot,\cdot\>_{t}$.
\end{Prop}

Therefore, the natural scalar product in the Kazhdan--Lusztig theory (the one for which the canonical basis is orthonormal) can be seen as a degenerate version of the Cherednik scalar product. One immediate consequence of these considerations is that parabolic Kazhdan--Lusztig polynomials can be obtained as
$$\<C'_\l,\widetilde E_\mu(\infty,t)\>_{t}=P^*_{v_\mu,v_\l}(t)$$
For root systems of type $A$ this is 
Lemma 11.3 in \cite{knop}. Similarly, from 
$$\<\widetilde E_\l(0,t),\widetilde E_\mu(\infty,t)\>_{t}=R^*_{v_\mu,v_\l}(t)$$
we obtain the parabolic $R$--polynomials (the notation is consistent with \cite{lusztig3}).
%%%%%%%%%%%%%%%%%%%%%%%%%%%%%%%%%%%%%%%%%%%%%%%%%%%%%%%%%%%%%%%%%%%%%%%%%%%%%
%%%%%%%%%%%%%%%%%%%%%%%%%%%%%%%%%%%%%%%%%%%%%%%%%%%%%%%%%%%%%%%%%%%%%%%%%%%%%
%%%%%%%%%%%%%%%%%%%%%%%%%%%%%%%%%%%%%%%%%%%%%%%%%%%%%%%%%%%%%%%%%%%%%%%%%%%%%

\section{Relating two limiting cases}\label{section5}

\subsection{The 0--Hecke algebra} The 0--Hecke algebra discussed here is a suitable degeneration of the Hecke algebra $\h$ as the parameters $t$ are specialized to zero.
\begin{Def}
The 0--Hecke algebra $\n$ associated to $\RR$ is the $\mathbb Q$--algebra described by generators and relations as follows:

\underline{Generators}: One generator $N_i$ for each simple root
$\a_i$.

\underline{Relations}: a) Each pair of generators satisfies the
same braid relations as the corresponding pair of simple reflections.

\hspace{1.52cm} b) The quadratic relations
$$
N_i^2 =  -N_i,  \quad 1\leq i\leq n
$$
\end{Def}

Since the generators $N_i$ satisfy the braid relations a standard basis of $\n$ is given by the elements $\{N_w\}_{w\in \W}$ where, as usual, $N_w=N_{i_l}\cdots N_{i_1}$ if 
$w=s_{i_l}\cdots s_{i_1}$ is a reduced  expression of $w$ in terms of simple reflections.
The 0--Hecke algebra has a linear action on $\cal R$ described by 
$$
N_i\cdot e^\l=-\frac{e^\l-e^{s_i(\l)}}{1-e^{-\a_i}}, \quad 1\leq i\leq n
$$

It is straightforward to check that the action of  $t_i^{\frac{1}{2}}T_i$  degenerates to the action of $N_i$ if we specialize the parameter $t_i$ to zero. In general $\chi(w)T_w$
will degenerate to $N_w$. In fact, by degenerating $t_i^{\frac{1}{2}}T_i^{-1}$ we obtain another set of generators $$N_i':=N_i+1 $$ satisfying the same braid relations and the quadratic relations $N_i'^2 =  N'_i$. Of course, $\chi(w)T^{-1}_{w^{-1}}$ will degenerate to $N'_w$. For this reason we call $\{N'_w\}_{w\in \W}$ the dual standard basis of $\n$. The operators $N_w$ are closely related to the Demazure operators $\Delta_w$, where $\Delta_i$ act as 
$$
\Delta_i\cdot e^\l=e^{s_i(\l)}+\frac{e^\l-e^{s_i(\l)}}{1-e^{-\a_i}}, \quad 1\leq i\leq n
$$
The relationship is
$$w_\circ N'_{w}w_\circ=\Delta_{w_\circ w w_\circ}$$

The generators $N'_i$ act on the standard basis as follows
\begin{eqnarray*}
N'_iN_w&=& N_w+N_{s_iw}, \quad \text{if}\ \ s_iw>w\\
N'_iN_w&=& 0,\quad \text{if}\ \ s_iw<w\\
\end{eqnarray*}
The following results are certainly well-known.
\begin{Lm} Let $w$ be an element of $\W$ and let $s_i$ be a simple reflection such that $s_iw>w$. Then,
$$
N'_i\sum_{x\le w} N_x=\sum_{y\le s_iw} N_y
$$

\end{Lm}
\begin{proof} Using the above two formulas we obtain
\begin{eqnarray*}
N'_i\sum_{x\le w} N_x &=& \sum_{x\le w, ~x<s_ix} (N_{s_ix}+N_x)\\
&=& \sum_{y\le s_iw} N_y
\end{eqnarray*}
The last equality followed from the third property of the Bruhat order.
\end{proof}
\begin{Cor}\label{cor10}
Let $w$ be an element of $\W$. Then,
$$
N'_w=\sum_{x\le w} N_x
$$
\end{Cor}
\begin{proof}
Apply the previous Lemma repeatedly.  
\end{proof}
%%%%%%%%%%%%%%%%%%%%%%%%%%%%%%%%%%%%%%%%%%%%%%%%%%%%%%%%%%%%%%%%%%%%%%%%%%%%%
%%%%%%%%%%%%%%%%%%%%%%%%%%%%%%%%%%%%%%%%%%%%%%%%%%%%%%%%%%%%%%%%%%%%%%%%%%%%%
%%%%%%%%%%%%%%%%%%%%%%%%%%%%%%%%%%%%%%%%%%%%%%%%%%%%%%%%%%%%%%%%%%%%%%%%%%%%%

\subsection{The limit $q\to 0$} In this section we collect some immediate consequences of Theorem \ref{state2} regarding the limit $q\to 0$.

\begin{Cor} Let $\l$ be a weight. Then, $$
\widetilde E_{\l}(0,t)=T_{\w_\l}\cdot e^{\l_-}
$$ In consequence, the coefficients of $E_\l(0,t)$ are polynomials in $\{t_\a\}_{\a\in R}$ with integer coefficients.
\end{Cor}
\begin{proof}
From Lemma \ref{lemma4} (2) we obtain that $X_{\l_-}=T_{w^{-1}_{\l_-}}^{-1}\omega_{\tilde \l}$. Hence, Theorem \ref{state2} implies 
$$
\widetilde E_{\l_-}(0,t)=e^{\l_-}
$$
From Lemma \ref{lemma4} (1) we deduce that $T^{-1}_{w_\l^{-1}}=T_{\w_\l}T^{-1}_{w^{-1}_{\l_-}}$. Now, Theorem \ref{state2} and the above formula give the desired result.
\end{proof}

\begin{Cor}\label{cor11} 
Let $\l$ be a weight. Then, $$
E_{\l}(0,0)=N_{\w_\l}\cdot e^{\l_-}
$$ 
\end{Cor}
\begin{proof}
The conclusion follows by sending the parameters $\{t_\a\}_{\a\in R}$ to zero in $$E_{\l}(0,t)=\chi(\w_\l) T_{\w_\l}\cdot e^{\l_-}$$
and keeping in mind that $\chi(\w_\l) T_{\w_\l}$ degenerate to $N_{\w_\l}$.
\end{proof}
Next, we explain the relationship between the $q\to 0$ limit and the $q\to \infty$ limit.

\begin{Prop} Let $\l$ be a weight. Then,
$$
w_\circ\cdot E_\l(\infty,t^{-1})=\chi(w_\circ\w_\l)T^{-1}_{(w_\circ\w_\l)^{-1}}\cdot e^{\l_-}
$$
\end{Prop}
\begin{proof} Let us argue first that $$
T^{-1}_{w_\circ}T^{-1}_{w^{-1}_\l}=T^{-1}_{(w_\circ \w_\l)^{-1}}T_{w_{\l_-}^{-1}}^{-1}
$$
Indeed, from Lemma \ref{lemma4} (1) we get $T_{w_\l}=T^{-1}_{\w_\l^{-1}}T_{w_{\l_-}}$. Keeping in mind that $T_{w_\circ}=T_{w_\circ \w_\l}T_{\w^{-1}_{\l}}$
we obtain $$
T_{w_\circ}T_{w_\l}=T_{w_\circ \w_\l}T_{w_{\l_-}}
$$
The claim  follows by applying the Kazhdan-Lusztig involution to this identity.

By Theorem \ref{state2},  $\widetilde E_\l(\infty,t)$ and $\widetilde E_\l(0,t)$ are interchanged by $\kappa$. This fact can be expressed as
$$
w_\circ\cdot\widetilde E_\l(\infty,t^{-1})=\chi(w_\circ)T^{-1}_{w_\circ}\cdot \widetilde E_\l(0,t)
$$
Now, 
\begin{eqnarray*}
T^{-1}_{w_\circ}\cdot \widetilde E_\l(0,t)&=&T^{-1}_{w_\circ}T^{-1}_{w^{-1}_\l}\cdot \widetilde E_{\tilde \l}(0,t)\\
&=& T^{-1}_{(w_\circ \w_\l)^{-1}}T_{w_{\l_-}^{-1}}^{-1}\cdot \widetilde E_{\tilde \l}(0,t)\\
&=& T^{-1}_{(w_\circ \w_\l)^{-1}}\cdot  e^{\l_-}
\end{eqnarray*}
The conclusion immediately follows. 
\end{proof}
Specializing further all the remaining parameters to 0 we obtain again the Demazure character formula.
\begin{Cor}
Let $\l$ be a weight. Then
$$
w_\circ\cdot E_\l(\infty,\infty)=N'_{w_\circ\w_\l}\cdot e^{\l_-}
$$
\end{Cor} 
To see that this is indeed equivalent to Demazure's formula 
$$
E_\l(\infty,\infty)=\Delta_{w_\l}\cdot e^{\l_+}
$$
note that $$w_\circ N'_{w_\circ \w_\l}w_\circ=\Delta_{\w_\l w_\circ}$$
The relationship between the limits $t\to0,\infty$ is described in the following
\begin{Cor}
Let $\l$ be a weight. Then
$$
w_\circ\cdot E_\l(\infty,\infty)=\sum_{w_\circ(\l)\le \mu\le \l_-} E_\mu(0,0)
$$
\end{Cor}
\begin{proof}
From Corollary \ref{cor10} we know that 
$$
N'_{w_\circ \w_\l}=\sum_{x\le w_\circ\w_\l} N_x
$$
Now, keep in mind that $N_i\cdot e^{\l_-}=0$ if $s_i$ fixes $\l_-$ and apply Corollary \ref{cor11} to obtain the desired result.
\end{proof}
%%%%%%%%%%%%%%%%%%%%%%%%%%%%%%%%%%%%%%%%%%%%%%%%%%%%%%%%%%%%%%%%%%%%%%%%%%%%%
%%%%%%%%%%%%%%%%%%%%%%%%%%%%%%%%%%%%%%%%%%%%%%%%%%%%%%%%%%%%%%%%%%%%%%%%%%%%%
%%%%%%%%%%%%%%%%%%%%%%%%%%%%%%%%%%%%%%%%%%%%%%%%%%%%%%%%%%%%%%%%%%%%%%%%%%%%%
\subsection{A geometric interpretation} \label{geometric} We first recall from \cite{ion} the geometric interpretation of the polynomials $E_\l(\infty,\infty)$.
For $w\in \W$ let $S_w$, respectively $S^-_w$, be the closure of the Bruhat cell $BwB/B$, respectively of $B^-wB/B$,  inside the flag variety $G/B$. For $\l$ a weight, $\cal L_\l$ 
denotes the corresponding line bundle over $G/B$. We will use the same notation for the restriction of this line bundle to any subvariety $S_w$ or $S_w^-$.

On the algebraic side, let $\l$ be a weight and let $V_{\l_+}$ the irreducible $G$--module with highest weight $\l_+$. By $V_{\l_+}(\l)$ we denote the (one dimensional) weight space of weight $\l$. The Demazure module corresponding to $\l$ is defined as 
$D_\l:=B\cdot V_{\l_+}(\l)$. The connection with geometry is the following: if $\l=w(\l_+)$, the Demazure module corresponding to $\l$ and the dual of the space of global sections of $\cal L_{-\l_+}$  over $S_w$ are isomorphic as $B$--modules
$$
H^0(S_w,\cal L_{-\l_+})^*\cong D_\l
$$
and the $T$--character of $D_\l$ is $E_\l(\infty,\infty)$. 

Equivalently, let $D^-_\l:=B^-\cdot V_{\l_+}(\l)$ and let $w$ such that $w(\l_+)=\l$. Then, 
$$
H^0(S^-_w,\cal L_{-\l_+})^*\cong D^-_\l
$$
as $B^-$--modules and the $T$--character of $D^-_\l$ is 
\begin{equation}\label{eq21}
w_\circ E_{w_\circ(\l)}(\infty,\infty)=\sum_{\l\le \mu\le \l_-} E_\mu(0,0)
\end{equation}

With the above notation, let $K_\l$ be the kernel of the restriction map 
$$
H^0(S^-_w,\cal L_{-\l_+})\to H^0(\bigcup_{w<y}S^-_y,\cal L_{-\l_+})
$$
(keep in mind that $S^-_y\subset S^-_w$ for $w\le y$). The dual of $K_\l$ is isomorphic to the co-kernel of the inclusion map 
$$
\bigcup_{\l<\mu\le \l_-} D^-_\mu \to D^-_\l
$$
By  (\ref{eq21}) the $T$-character of $\bigcup_{\l<\mu\le \l_-} D^-_\mu$ equals $$\sum_{\l< \mu\le \l_-} E_\mu(0,0)$$ and therefore
the character of $K_\l^*$ is $E_\l(0,0)$.
\begin{Thm}\label{00}
Let $w$ be an element of $\W$ and let $\l=w(\l_+)$. Then, $E_\l(0,0)$ is the character of the dual space of sections of $\cal L_{-\l_+}$ which are supported on $$S^-_w-\bigcup_{w<y}S^-_y$$
In consequence, they are polynomials with non-negative integer coefficients.
\end{Thm}

As a terminological coincidence, in type $A$ the earliest reference to the polynomials $E_\l(0,0)$ seems to go back to the work of Lascoux and Sch\" utzenberger  \cite[Theorem 3.8]{lascoux} where they form their ``standard basis".
%%%%%%%%%%%%%%%%%%%%%%%%%%%%%%%%%%%%%%%%%%%%%%%%%%%%%%%%%%%%%%%%%%%%%%%%%%%%%
%%%%%%%%%%%%%%%%%%%%%%%%%%%%%%%%%%%%%%%%%%%%%%%%%%%%%%%%%%%%%%%%%%%%%%%%%%%%%
%%%%%%%%%%%%%%%%%%%%%%%%%%%%%%%%%%%%%%%%%%%%%%%%%%%%%%%%%%%%%%%%%%%%%%%%%%%%%

\end{document}